\documentclass[pdflatex,sn-mathphys-ay]{sn-jnl}% Math and Physical Sciences Author Year Reference Style
\usepackage{graphicx}%
\usepackage{multirow}%
\usepackage{amsmath,amssymb,amsfonts}%
\usepackage{amsthm}%
\usepackage{mathrsfs}%
\usepackage[title]{appendix}%
\usepackage{xcolor}%
\usepackage{textcomp}%
\usepackage{manyfoot}%
\usepackage{booktabs}%
\usepackage{algorithm}%
\usepackage{algorithmicx}%
\usepackage{algpseudocode}%
\usepackage{listings}%

\usepackage[utf8]{inputenc}
\usepackage{dsfont}
\theoremstyle{thmstyleone}%
\newtheorem{theorem}{Theorem}%  meant for continuous numbers
\newtheorem{proposition}[theorem]{Proposition}% 
\newtheorem{assumption}{Assumption}

\theoremstyle{thmstyletwo}%

\theoremstyle{thmstylethree}%

\newcommand{\eps}{{\varepsilon}}
\renewcommand{\phi}{\varphi}
\newcommand{\R}{\mathbb{R}}

\newcommand{\Z}{\mathbb{Z}}
\newcommand{\N}{\mathbb{N}}

\newcommand{\pr}{\mathbb{P}}       
\newcommand{\ex}{\mathbb{E}}       
\newcommand{\var}{\textnormal{Var}} 
\newcommand{\cov}{\textnormal{Cov}}
\newcommand{\Nc}{\mathcal{N}}
\newcommand{\Uc}{\mathcal{U}}
\newcommand{\Fc}{\mathcal{F}}

\newcommand{\Gc}{\mathcal{G}}

\newcommand{\Oc}{\mathcal{O}}

\newcommand{\Sb}{\mathbf{S}}

\newcommand{\diff}{{\,\mathrm{d}}}
\newcommand{\convw}{\rightsquigarrow}                           
 
\newcommand{\id}{\mathds{1}}

\DeclareMathOperator*{\argmin}{argmin}

\begin{document}
	
	\title[Functional CLT for Localized Partial Sums]{A Functional Central Limit Theorem for Localized Partial Sums of Non-Stationary Time Series}	
	
	%\author*[1,2]{\fnm{Florian} \sur{Heinrichs}}\email{f.heinrichs@fh-aachen.de}
	\author{\fnm{Florian} \sur{Heinrichs}}\email{f.heinrichs@fh-aachen.de}
	
	\affil{\orgdiv{Department of Medical Engineering and Technomathematics}, \orgname{FH Aachen - University of Applied Sciences}, \orgaddress{\street{Heinrich-Mußmann-Straße 1}, \city{Jülich}, \postcode{52428}, \country{Germany}}}
	
	%\affil[2]{\orgdiv{Department}, \orgname{Organization}, \orgaddress{\street{Street}, \city{City}, \postcode{10587}, \state{State}, \country{Country}}}
	
	%\affil[3]{\orgdiv{Department}, \orgname{Organization}, \orgaddress{\street{Street}, \city{City}, \postcode{610101}, \state{State}, \country{Country}}}
		
	\abstract{A localized functional central limit theorem is established for kernel-weighted partial sum processes of piecewise locally stationary time series under geometric decay of the physical dependence measure. 
	The localized process is shown to converge weakly to a centered Gaussian random distribution in $D'(0,1)$, and the limit extends naturally to an isonormal Gaussian process on $L^2([0,1])$.
	Weak convergence is further derived for processes indexed by totally bounded subsets of $L^2([0,1])$.
	As an application, the localized limit theory is used to construct tests for constant mean functions against linear, polynomial, and general alternatives in non-parametric regression with locally stationary errors.
	Simulation results and data examples illustrate the finite sample performance and practical applicability of the proposed methodology.}
	
	\keywords{Physical Dependence Measure, Piecewise Local Stationarity, Functional Central Limit Theorem; Kernel-Weighted Partial Sums; Change Detection}
	
	\pacs[MSC Classification]{60F17, 62M10, 60G65, 62G20}
	
	\maketitle
	
	% !TeX spellcheck = en_US
%% ====================
\section{Introduction} \label{sec:intro}
%% ====================

% Kernel regression -> Convolution of time series (e.g. mean estimation)
% Formula: $\sum_{i=1}^n ...$
% Pointwise CLT (cite Dahlhaus and others)
% Empirical process (cite Richter)
% For separate points: Limits indepedent (white noise) -> no FCLT should be expected
% But: There is more structure
%  -> For points close together, dependence
% Two solutions:
%  - either "local" CLT at h_n scale: H_n(t, s) = G_n(t + s h_n), but this leads to independence among t, so again (functional) white noise
%  - alternatively: work with ``distributions'' and generalized Gaussian processes

% Statistical applications: 
%  - Testing for changes in \mu (random projections)
%  - Goodness-of-fit in parametric model \mu(x) = a x + b

%\begin{equation*}
%	G_n(t) = \frac{1}{\sqrt{n}h_n} \sum_{i=1}^{n} (X_{i, n} - \ex[X_{i, n}]) K_{h_n}(\tfrac{i}{n}-t) \convw G(t),
%\end{equation*}
%for some Gaussian process $G\in D'(0, 1)$.

Modern time series data are frequently non-stationary, with serial dependence structures that evolve over time through gradual or abrupt changes in distributional characteristics such as the mean, variance, or autocovariance.  
Examples arise naturally in climatology, neuroscience, economics, and engineering, where a global stationarity assumption is often too restrictive, while a fully unrestricted non-stationary model is too broad to support useful asymptotic theory.  
A particularly successful compromise is the framework of local stationarity and its extensions, which approximate a non-stationary process locally by a family of stationary processes and thereby permit a rich asymptotic analysis of estimators and test statistics.

Within this broad area, the physical dependence measure introduced by \cite{wu2005} has become one of the most powerful tools for deriving limit theorems for weakly dependent processes generated as non-linear functionals of an underlying i.i.d. innovation sequence.  
Unlike purely covariance-based notions of dependence, the physical dependence measure is tied directly to the data generating mechanism and is therefore well suited to non-linear and non-stationary settings.  
Its flexibility has led to a large body of work on Gaussian approximation, invariance principles, and statistical inference for both stationary and locally stationary time series, including the piecewise locally stationary framework considered in this paper.  
For the present results, this perspective is particularly natural because our main assumptions and proofs are built explicitly on physical dependence bounds and on the Gaussian approximation theory derived from them. 

We study the localized partial sum process  
\begin{equation*}
	G_n(t) = \frac{1}{\sqrt{n}h_n} \sum_{i=1}^{n} (X_{i, n} - \ex[X_{i, n}]) K_{h_n}(\tfrac{i}{n}-t), \qquad t\in[0, 1],
\end{equation*}
where $K$ is a compactly supported kernel and $h_n\to 0$ is a bandwidth sequence.  
Processes of this form arise naturally in nonparametric smoothing, local change detection, and kernel-based estimation in non-stationary environments.  
At a fixed rescaled time point $t$, the statistic $G_n(t)$ aggregates observations from a shrinking neighborhood of $t$, so its asymptotic behavior differs fundamentally from that of global partial sums.  
In particular, localized sums are sensitive to the local long-run variance structure and can be used to probe local deviations of the underlying mean function after suitable smoothing and bias correction.

Pointwise asymptotic normality for such localized quantities is by now reasonably well understood under local stationarity.  
In particular, a pointwise central limit theorem for localized kernel-weighted sums follows from existing theory for non-linear locally stationary processes, see for example \cite{dahlhaus2019}.  
However, pointwise weak convergence does not by itself yield a satisfactory process-level description of the entire stochastic process $(G_n(t))_{t\in[0, 1]}$.  
For many statistical applications, including integrated or projected test statistics and inference over classes of functions, one needs a functional central limit theorem that captures the joint asymptotics of all linear functionals of the localized process.

The main contribution of this paper is such a localized functional central limit theorem for piecewise locally stationary time series in the sense of \cite{zhou2013}.  
More precisely, we establish weak convergence of $G_n$ in the space of distributions $D'(0,1)$, viewed through its action on test functions.  
The limit is a centered Gaussian random distribution whose covariance structure is determined by the time-varying long-run variance function $\sigma^2(\cdot)$.  
At continuity points of $\sigma^2$, the result is consistent with the usual local variance formula from pointwise theory, while the distribution-valued formulation allows us to move beyond pointwise asymptotics and to treat integrated and projected versions of $G_n$.

Our proof strategy relies on a Gaussian approximation for piecewise locally stationary processes established by \cite{zhou2013}, which extends the earlier Gaussian approximation by \cite{wu2011}. The former is based on physical dependence arguments in the spirit of \cite{zhou2009}.  
This approximation allows us to replace the localized partial sums by kernel-weighted Gaussian sums with time-varying variance, after which the limiting covariance structure becomes explicit.  
The resulting theorem shows that physical dependence methods are strong enough not only for pointwise normal approximation, but also for a distribution-valued functional limit theory in a non-stationary setting.

A second contribution is an extension from test-function convergence to convergence of finite-dimensional distributions indexed by $L^2([0,1])$-functions.  
This yields a natural isonormal Gaussian limit representation, where the time-varying long-run variance enters through a multiplication operator.  
We further provide a sufficient condition for asymptotic tightness of the induced process on totally bounded subsets of $L^2([0,1])$, thereby obtaining weak convergence in $\ell^\infty(\Gc)$ for suitable function classes $\Gc$.  
This bridges the distribution-valued weak limit with empirical-process style statements indexed by structured families of test functions.

The present work is related to several strands of literature.  
First, it contributes to the theory of locally stationary and piecewise locally stationary time series initiated in seminal work by \cite{dahlhaus1996} and further developed in many directions \citep{vogt2012,zhou2009}. For the specific PLS framework used here we rely on the formulation by \cite{zhou2013} and dependence conditions based on the physical dependence measure in the spirit of \cite{wu2005}.
Second, our point of departure is the pointwise CLT for localized statistics under non-linear local stationarity, see e.\,g. \cite{dahlhaus2019}, which we lift to a process-level statement in the space of distributions.  
Third, our results are connected to recent work on empirical processes for locally stationary time series. \cite{phandoidaen2022} studied weak convergence of the generalized empirical process of multivariate locally stationary time series under physical dependence. More recently, \cite{palm2025} provided similar results for locally stationary time series in Polish spaces under $\beta$-mixing.
These works address complementary questions on process convergence and empirical process methodology under local stationarity, whereas our focus is on kernel-localized partial sums and their Gaussian distribution-valued limit.  
A further related line concerns change point and trend detection for non-stationary time series, including bootstrap and self-normalization approaches under physical dependence conditions \citep{zhou2013,bucher2021,heinrichs2021,heinrichs2025b}.
Finally, the localized partial sum process $G_n$ naturally occurs in a variety of contexts, such as outlier detection \citep{heinrichs2025}.

Beyond the theoretical contribution, the localized functional central limit theorem developed here has direct statistical consequences.  
We apply it to inference in the nonparametric regression model  
\begin{equation*}
	X_{i,n}=\mu(i/n)+\eps_{i,n},
\end{equation*}
where $\mu$ is an unknown mean function and $(\eps_{i,n})$ are centered locally stationary errors.  
Using bias-reduced local linear smoothing, we derive asymptotically valid tests for constancy of the mean against linear, polynomial, and general alternatives.  
The resulting procedures are formulated through projections of the estimated mean onto deterministic or random test functions and are justified by the Gaussian limit theory for localized partial sums.

Our simulation study compares the proposed methods with the bootstrap-based structural change test of \cite{zhou2013} and with a CUSUM benchmark.  
The results indicate that the new procedures can be competitive and often powerful, although they display some finite-sample conservatism or size distortion depending on the form of non-stationarity and serial dependence.  
We further illustrate the methodology with applications to daily temperature data and EEG recordings, highlighting the practical relevance of localized testing procedures under non-stationarity.

The paper is organized as follows.  
Section \ref{sec:preliminaries} introduces the framework of piecewise local stationarity, the physical dependence measure, and the distribution-space setting used for the functional limit theorem.  
Section \ref{sec:method} contains the main theoretical results, including the localized functional central limit theorem in $D'(0,1)$, the $L^2$-indexed Gaussian extension, and a tightness criterion for function-indexed processes.  
Section \ref{sec:stat_app} develops applications to testing for changes in the mean function under locally stationary errors.  
Section \ref{sec:empirical} reports simulation results and data examples, and Section \ref{sec:proofs} contains the proofs.

	% !TeX spellcheck = en_US
%% ====================
% \section{Main Result} \label{sec:method}
%% ====================

\section{Mathematical Preliminaries} \label{sec:preliminaries}

\textbf{Random Distributions.} Recall that a sequence of random variables $(X_n)_{n\in\N}$ converges weakly to $X$ in some metric space $E$, whenever $\lim_{n\to\infty}\ex[f(X_n)] = \ex[f(X)]$ for any bounded, continuous function $f:E \to \R$. For random variables in $\R$, weak convergence implies tightness. Moreover, weak convergence is equivalent to pointwise convergence of the characteristic functions
\begin{equation*}
	\phi_n(t) = \ex[\exp(it X_n)] \xrightarrow{n\to\infty} \ex[\exp(it X)] = \phi(t),
\end{equation*}
for every $t\in\R$, by Lévy's continuity theorem. For stochastic processes in $\ell^\infty([0, 1])$, weak convergence is equivalent to weak convergence of marginal distributions and tightness of the process \citep[see, e.\,g., Theorem 1.5.4 in][]{vandervaart2023}. 

Let $C_c^\infty(0, 1)$ denote the space of all infinitely differentiable functions on $(0, 1)$ with compact support, and let $D'(0, 1)$ denote the strong dual space of $C_c^\infty(0,1)$. With this notation, $C_c^\infty(0, 1)$ is referred to as the \textit{space of test functions on $(0, 1)$} and $D'(0, 1)$ as the \textit{space of distributions on $(0, 1)$}. 

Let $\langle f, g\rangle = \int_0^1 f(x) g(x) \diff x$ and let ``$\convw$'' denote weak convergence. A sequence of ``stochastic processes'' $(X_n)_{n\in\N}$ in $D'(0, 1)$ converges weakly to $X \in D'(0, 1)$, iff $\langle X_n, \phi \rangle \convw \langle X, \phi \rangle$ for any test function $\phi \in C_c^\infty(0, 1)$. This definition of weak convergence seems substantially weaker than usual weak convergence in metric spaces. However, it satisfies the key properties discussed before.

First, weak convergence in $D'(0, 1)$ is equivalent to pointwise convergence of characteristic functionals, i.\,e.,
\begin{equation*}
	\ex[\exp(i \langle X_n, \phi \rangle)] \xrightarrow{n\to\infty} 	\ex[\exp(i \langle X, \phi \rangle)],
\end{equation*}
for any $\phi \in C_c^\infty(0, 1)$, by Theorem III.6.5 of \cite{fernique1967}. The latter theorem is the generalized distribution-analog of Lévy's continuity theorem. 
Moreover, by bilinearity of $\langle \cdot, \cdot\rangle$ and the Cramér-Wold device, weak convergence in $D'(0, 1)$ is equivalent to weak convergence of the marginal distributions, i.\,e.,
\begin{equation*}
	( \langle \phi_1, X_n \rangle, \dots, \langle \phi_d, X_n \rangle )
	\convw ( \langle \phi_1, X \rangle, \dots, \langle \phi_d, X \rangle ).
\end{equation*} 
%
% Finally, by Mitoma's theorem \citep[see, e.\,g., Theorem 6.13 in ][]{walsh1986}, tightness in $D'(0, 1)$ follows by considering the constant process $s\mapsto X_n$ in $D([0, 1], D'(0, 1))$.
Finally, applying Mitoma's theorem \citep[see, e.,g., Theorem 6.13 in ][]{walsh1986} to the constant processes $\tilde{X}_n(s) =X_n$, $s\in[0,1]$, yields tightness of $(\tilde{X}_n)_{n\in\N}$ in $D([0,1],D'(0,1))$ whenever $(\langle X_n,\phi\rangle)_{n\in\N}$ is tight in $\R$ for every $\phi\in C_c^\infty(0,1)$. Since the constant-path embedding is a topological embedding of $D'(0,1)$, equipped with its strong topology, into $D([0,1],D'(0,1))$ with the Skorokhod topology, this implies tightness of $(X_n)_{n\in\N}$ in $D'(0,1)$.

Therefore, $\langle X_n, \phi \rangle \convw \langle X, \phi \rangle$ for any test function $\phi \in C_c^\infty(0, 1)$, is indeed the suitable generalization of weak convergence to random elements of $D'(0, 1)$. 

Similarly, the distribution of an element $G \in D'(0, 1)$ is uniquely determined by its marginals $( \langle \phi_1, G \rangle, \dots, \langle \phi_d, G \rangle )$, for any test functions $\phi_1, \dots, \phi_d \in C_c^\infty(0, 1)$ and $d\in\N$ \citep[see, e.\,g., Proposition III.4.2(b) of ][]{fernique1967}.

\textbf{Piecewise Local Stationarity.} If the time series $X = \{(X_{i,n})_{i=1, \dots, n}\}_{n\in \N}$ is arbitrarily irregular, we should not expect a central limit theorem to hold. In the following, we work within the framework of piecewise local stationarity as introduced by \cite{zhou2013}.

Let $(\eta_i)_{i\in\Z}$ denote a sequence of i.i.d. random variables and let $\Fc_i = (\eta_k)_{k \le i}$. A time series is \textit{piecewise locally stationary} (PLS) with $r$ break points if there exist time points $0 = t_0 < t_1 < \dots < t_r < t_{r+1} = 1$ and (possibly nonlinear) filters $H_0, H_1, \dots, H_r$, with $H_j:[0, 1] \times \R^\N \to \R$, for $j=0,\dots, r$, such that
\begin{equation*}
	X_{i, n} = H_j(\tfrac{i}{n}, \Fc_i), \quad \mathrm{if}~t_j < \tfrac{i}{n} \le t_{j+1}.
\end{equation*}

A PLS time series with $0$ break points is locally stationary in the sense of \cite{zhou2009}. If $H_j(t, \cdot) \equiv H_j(0, \cdot)$ for all $j=0, \dots, r$ and $t\in[0,1]$, then the time series is piecewise stationary. In particular, PLS allows for abrupt and gradual changes in the time series' distribution. 

A PLS time series is \textit{piecewise Lipschitz continuous} in $L^2$-norm, if for all $j= 0, \dots, r$, and $t, s \in [t_j, t_{j+1}]$,
\begin{equation*}
	\ex[|H_j(t, \Fc_0) - H_j(s, \Fc_0)|^2]^{1/2} \le C |t - s|,
\end{equation*}
for some constant $C\ge 0$.

A natural measure of dependence of PLS time series, is the \textit{physical dependence measure}, as originally introduced by \cite{wu2005}, and extended to PLS by \cite{zhou2013}. Let $(\eta_i^*)_{i\in\N}$ denote an independent copy of $(\eta_i)_{i\in\N}$, and let $\Fc_i^* = (\dots, \eta_{-1}, \eta_0^*, \eta_1, \dots, \eta_i)$ be the filtration with $\eta_0$ replaced by its independent copy $\eta_0^*$. Then, the physical dependence measure of order $p$ at lag $h \ge 0$ is defined by
\begin{equation*}
	\delta_p(h) = \max_{0\le j \le r} \sup_{t\in[t_j, t_{j+1}]} \ex[|H_j(t, \Fc_h) - H_j(t, \Fc_h^*)|^p]^{1/p}.
\end{equation*}

\section{Localized Functional Central Limit Theorems} \label{sec:method}

When considering the \textit{localized partial sum process} 
\begin{equation*}
	G_n(t) = \frac{1}{\sqrt{n}h_n} \sum_{i=1}^{n} (X_{i, n} - \ex[X_{i, n}]) K_{h_n}(\tfrac{i}{n}-t),
\end{equation*}
for $t\in[0, 1]$ and a suitable kernel $K$, we must distinguish between different asymptotic regimes. While pointwise weak convergence of $\sqrt{h_n}G_n(t)$ in $\R$ to a normal distribution is well understood \citep[see, e.\,g., Theorem 2.10 in][]{dahlhaus2019}, we establish uniform convergence of $G_n$ in $D'(0, 1)$ under the following regularity assumptions. 

\begin{assumption} \label{assump:kern}
	The kernel $K:\R\to\R$ is symmetric, supported on $[-1, 1]$ and Lipschitz continuous on its support with $\int K(x) \diff x = 1$.
\end{assumption}

\begin{assumption} \label{assump:error}
	The piecewise locally stationary time series $X = \{(X_{i,n})_{i=1, \dots, n}\}_{n\in \N}$ satisfies the following conditions.
	\begin{enumerate}
		\item $X$ is piecewise Lipschitz continuous in $L^2$-norm.
		\item $X$ has finite fourth moments, i.\,e., $\ex[|H_j(t, \Fc_0)|^4]^{1/4} < \infty$ for all $j=0,\dots, r$ and $t\in[t_j, t_{j+1}]$.
		\item $\delta_4(h) = \Oc(\chi^h)$ for some $\chi\in(0, 1)$.
		\item The long-run variance of $X$, defined by
		\begin{equation*}
			\sigma^2(t) = \sum_{k\in\Z} \cov\big( H_j(t, \Fc_0), H_j(t, \Fc_k) \big),
		\end{equation*}
		for $t\in[t_j, t_{j+1}]$, is bounded away from $0$, i.\,e., $\inf_{t\in[0, 1]}\sigma^2(t) > 0$.
	\end{enumerate}
\end{assumption}

% \begin{remark}
	Assumption \ref{assump:kern} is rather a design choice for the localized partial sum process and satisfied by typical kernels, such as the uniform, triangular, Epanechnikov and quartic kernels.
	Assumption \ref{assump:error} is typical under the framework of piecewise local stationarity \citep[see, e.\,g.,][]{zhou2013}. 
% \end{remark}

\begin{theorem} \label{thm:main}
	Let $(h_n)_{n\in\N}$ be a sequence with $\lim_{n\to\infty} h_n = 0$. Under Assumptions \ref{assump:kern} and \ref{assump:error}, the following convergences hold:
	\begin{enumerate}
		\item If $n^{1/4}h_n^{1/2}/\log^2n \to \infty$, $\sqrt{h_n} G_n(t) \convw \Nc(0, \sigma_K^2(t))$, for each $t\in (0, 1)$,
		where 
		\begin{equation*}
			\sigma_K^2(t) = \sigma^2(t_-) \int_{-1}^0 K^2(x) \diff x + \sigma^2(t_+) \int_0^1 K^2(x) \diff x.
		\end{equation*}
		In particular, $\sigma_K^2(t) = \sigma^2(t) \int K^2(x)\diff x$ for continuity points of $\sigma^2$.		
		\item If $n^{1/4}h_n/\log^2n \to \infty$, $G_n \convw G$ in $D'(0, 1)$, where $G$ is a centered Gaussian random distribution characterized by its finite-dimensional distributions 
		\begin{equation*}
			( \langle \phi_1, G \rangle, \dots, \langle \phi_d, G \rangle ) \sim \Nc(0, \Sigma),
		\end{equation*}
		where $\Sigma = (\Sigma_{i,j})_{i, j=1}^d$ is defined by $\Sigma_{i,j} = \int_0^1 \sigma^2(x) \phi_i(x) \phi_j(x) \diff x$.
		% \begin{equation*}
		% 	\Sigma_{i,j} = \int_0^1 \sigma^2(x) \phi_i(x) \phi_j(x) \diff x.
		% \end{equation*}
	\end{enumerate}
	 
\end{theorem}

Assumption \ref{assump:error} is required to apply Proposition 5 of \cite{zhou2013}, which extends the Gaussian approximation from Corollary 1 of \cite{wu2011} to PLS time series. The latter only requires moments of order $2+\eps$ and a geometric decay of $\delta_{2+\eps}(h)$, for some $\eps \in (0, 2]$. As long as $n^{\eps/(4+2\eps)}h_n/\log^2n \to \infty$, the theorem remains valid by an analogous Gaussian approximation.

The pointwise statement in part 1 is restricted to interior points $t\in(0,1)$, since kernel localization induces the usual boundary effects near 0 and 1. By contrast, the distribution-valued convergence in part 2 is formulated through test functions $\phi\in C_c^\infty(0,1)$, whose compact support excludes boundary contributions. Theorem 2 then extends this result to $L^2([0,1])$ by density, and Theorem 3 builds on that extension. 

Recall that $W = \{W(f)\}_{f\in L^2([0, 1])}$ is an isonormal Gaussian process on $L^2([0, 1])$, if $\ex[W(f)] = 0$ and $\ex[W(f)W(g)] = \int f(x) g(x)\diff x$, for any $f,g \in L^2([0, 1])$ \citep{dudley1973}. The limit $G \in D'(0, 1)$ can be extended to such an isonormal Gaussian process. Intuitively, $C_c^\infty(0, 1)$ is dense in $L^2([0, 1])$, and we may define $W(f) = L^2(\Omega)-\lim_{\phi\to f} \langle \phi, G\rangle$, for $f \in L^2([0, 1])$. Then, convergence of the finite-dimensional distributions of $G_n$ to $G$ extends to the corresponding convergence to $W$.

This extension is crucial for the statistical applications developed later. In Section \ref{sec:stat_app}, the test statistics are formulated as projections of the localized process onto deterministic or random functions in $L^2([0,1])$, rather than onto smooth test functions in $C_c^\infty(0,1)$. Theorem \ref{thm:isonormal} therefore provides the asymptotic justification for these projected statistics and links the distribution-valued limit in $D'(0,1)$ to the projection-based testing framework.

\begin{theorem} \label{thm:isonormal}
	Under the assumptions of Theorem \ref{thm:main} (2), for any $\phi_1, \dots, \phi_d \in L^2([0, 1])$,
	\begin{equation*}
		(\langle \phi_1, G_n \rangle, \dots, \langle \phi_d, G_n \rangle) \convw (W_\sigma(\phi_1), \dots, W_\sigma(\phi_d)),
	\end{equation*}
	where $W_\sigma = W \circ T_\sigma$, $T_\sigma(f) = f \cdot \sigma$ and $W = \{W(f)\}_{f\in L^2([0, 1])}$ denotes an isonormal Gaussian process on $L^2([0, 1])$. 
\end{theorem}

% \begin{remark}
	Equivalently, $W_\sigma$ is an isonormal Gaussian process on $L^2([0, 1], d_\sigma)$, where $d_\sigma$ denotes the finite Borel measure on $[0, 1]$ with density $\sigma^2$ with respect to the Lebesgue measure.
% \end{remark}

% Let $\Gc = C^\alpha([0, 1])$ denote the Hölder space for some $\alpha > 0$, equipped with the norm $\|f\|_\alpha = \|f\|_\infty + \sup_{s\neq t} \frac{|f(s) - f(t)|}{|s-t|}$. Define the process $H_n: \Gc \to \R$ by $H_n(f) = \langle G_n, f\rangle$. 

As discussed earlier, weak convergence of a stochastic process is equivalent to convergence of the finite dimensional distributions and (asymptotic) tightness of the process. Let $H_n =\{H_n(f)\}_{f\in\Gc}$ be the process defined by $H_n(f) = \langle f, G_n \rangle$ for some family of functions $\Gc \subset L^2([0, 1])$. If $\Gc$ is large, e.\,g., $\Gc = L^2([0, 1])$, we should not expect $H_n$ to be tight, so that Theorem \ref{thm:isonormal} does not extend to weak convergence of $H_n$.

While Theorem \ref{thm:isonormal} yields convergence of arbitrary finite collections of $L^2([0,1])$-projections, many statistical procedures involve a whole family of projection functions rather than finitely many fixed directions. This is relevant, for example, when one wishes to optimize over a class of test functions or to study statistics indexed by subsets of $L^2([0,1])$. Theorem \ref{thm:hoelder} provides a sufficient condition on $\Gc$ for asymptotic tightness of $H_n$, and hence, weak convergence of $H_n$ in $\ell^\infty(\Gc)$ to a Gaussian limit, thereby complementing the finite-dimensional extension from Theorem \ref{thm:isonormal} by a genuine process-level result.

\begin{theorem} \label{thm:hoelder}
	Let $\Gc$ be a totally bounded subset of $L^2([0, 1])$, $(h_n)_{n\in\N}$ a sequence with $\lim_{n\to\infty} h_n = 0$ and Assumptions \ref{assump:kern} and \ref{assump:error} be satisfied. Let $\{H_n(f)\}_{f\in\Gc}$ be defined by $H_n(f) = \langle f, G_n \rangle$.
	If $n^{1/4}h_n/\log^2n \to \infty$, $H_n \convw H$ in $\ell^\infty(\Gc)$, where $H$ is a centered Gaussian process with covariance
	\begin{equation*}
		\cov(H(f), H(g)) = \int_0^1 \sigma^2(x) f(x)g(x) \diff x.
	\end{equation*}
\end{theorem}

% \begin{remark}
	The unit ball in the Hölder space $C^\alpha([0, 1])$ is totally bounded with respect to the $L^2([0, 1])$-norm, such that the theorem holds. More specifically, the unit ball is bounded and equicontinuous, hence its closure is compact in $C([0, 1])$ by Arzelà-Ascoli. The embedding $C([0,1])\hookrightarrow L^2([0,1])$ is continuous, such that the image of the unit ball of $C^\alpha([0, 1])$ is compact in $L^2([0, 1])$, hence, totally bounded.
% \end{remark}

\section{Application to Change Detection} \label{sec:stat_app}
Localized partial sum processes appear naturally in the context of kernel regression. Consider the additive model
\begin{equation*}
	X_{i, n} = \mu(i/n) + \eps_{i, n},\quad i= 1,\dots, n,
\end{equation*}
for a smooth function $\mu$ and a centered, locally stationary error process $\eps$. Then, we are often interested in testing
\begin{equation} \label{eq:hypotheses}
	H_0: \mu(t) \equiv \mu(0)~\mathrm{for~all}~t\in[0, 1] \quad \mathrm{vs.} \quad H_1: \mu(t) \neq \mu(0)~\mathrm{for~some}~t\in[0, 1].
\end{equation}
In order to test whether the unknown mean $\mu$ is constant, we must first estimate it. Let $K$ be a kernel and $(h_n)_{n\in\N}$ a positive bandwidth with $h_n\to 0$, and let $K_{h_n}(\cdot)=K(\cdot/h_n)$. The local linear estimator of $\mu$ with bandwidth $h_n$ is defined as the first coordinate of
\begin{equation*}
	\big(\hat{\mu}_{h_n}(t), \widehat{\mu'}_{h_n}(t)\big) = \argmin_{b_0, b_1\in\R}\sum_{i=1}^n \bigg(X_{i,n} - b_0 - b_1 \big(\tfrac{i}{n} - t\big)\bigg)^2 K_{h_n}\big(\tfrac{i}{n} - t\big),
\end{equation*}  
see, for example, \cite{fan1996}. The local linear estimator of $\mu$ has a bias term proportional to $h_n^2$. Using the Jackknife bias reduction technique, as proposed by \cite{schucany1977}, we define $\tilde{\mu}_n(t) = 2 \hat{\mu}_{h_n/\sqrt{2}}(t) - \hat{\mu}_{h_n}(t)$, 
% \begin{equation*}
% 	\tilde{\mu}_n(t) = 2 \hat{\mu}_{h_n/\sqrt{2}}(t) - \hat{\mu}_{h_n}(t),
% \end{equation*}
to remove the leading order bias. Then $K^*(x) = 2 \sqrt{2}K(\sqrt{2}x) - K(x)$ is the kernel corresponding to $\tilde{\mu}_n$. If $\mu$ is twice differentiable with Lipschitz continuous second derivative, 
\begin{equation*}
	\sup_{t\in[h_n, 1-h_n]} | \tilde{\mu}_n(t) - \mu(t) - n^{-1/2}G_n^\eps(t)| = O(h_n^3 + \tfrac{1}{nh_n}),
\end{equation*}
where $G_n^\eps(t) = (\sqrt{n}h_n)^{-1}\sum_{i=1}^{n} \eps_{i, n} K^*_{h_n}(i/n - t)$ denotes the localized partial sum process of $\eps$ w.r.t. $K^*$, by Lemma C.2 in the Supplementary Material of \cite{dette2019}.

\subsection{Testing for a Linear Trend} \label{sec:lin_test}

In the simplest case, we assume $\mu(t) = \alpha_0 + \alpha_1 t$, such that the testing problem is equivalent to
\begin{equation}\label{eq:hypotheses_lin}
	H_0^{(\mathrm{lin})}: \alpha_1 = 0 \quad \mathrm{vs.} \quad H_1^{(\mathrm{lin})}: \alpha_1 \neq 0.
\end{equation}
Let $\phi_1(t) = t - 1/2$ and define the test statistic
\begin{equation} \label{eq:stat_lin}
	S_n(\phi_1) = \sqrt{n} \langle \phi_1, \tilde{\mu}_{n}\rangle. % \int_0^1 \phi_1(t) \tilde{\mu}_{n}(t) \diff t.
\end{equation}
Then, since $\int_0^1 \phi_1(t) \diff t = 0$ by definition, 
\begin{equation} \label{eq:lin_expansion}
	% S_n(\phi_1) = \sqrt{n} \int_0^1 \phi_1(t) \big(\tilde{\mu}_{n}(t) - \mu(t)\big) \diff t + \sqrt{n} \int_0^1 \phi_1(t) \bigg(\mu(t) - \int_0^1 \mu(x)\diff x\bigg) \diff t.
	S_n(\phi_1) = \sqrt{n} \langle \phi_1,  \tilde{\mu}_{n} - \mu \rangle + \sqrt{n} \Big\langle \phi_1, \mu - \int_0^1 \mu(x)\diff x\Big\rangle.
\end{equation}
The first term on the right-hand side converges to a centered normal distribution with variance $\int_0^1 \sigma^2(t) \phi_1^2(t) \diff t$. For the second summand,
\begin{align} \notag
	\Big\langle \phi_1, \mu - \int_0^1 \mu(x)\diff x\Big\rangle
	& = \int_0^1 \phi_1(t) \bigg(\alpha_0 + \alpha_1 t - \int_0^1 \alpha_0 + \alpha_1 x \diff x\bigg) \diff t \\
	& = \alpha_1 \int_0^1 \phi_1(t) (t - \tfrac{1}{2}) \diff t, \label{eq:lin_expansion2}
\end{align}
which equals zero under $H_0$ and is different from $0$ under the alternative. 

To derive critical values of $S_n(\phi_1)$ under the null hypothesis, we need an estimator of $\sigma^2(t)$. Define the weight function $w_{\tau_n}(t, i) = K_{\tau_n}(i/n - t) / \{\sum_{j=1}  K_{\tau_n}(j/n - t)\}$, for some positive bandwidth $\tau_n = o(1)$, and let $(m_n)_{n\in\N}$ be an integer-valued sequence with $m_n \to \infty$ and $m_n = o(n)$. Define
\begin{equation} \label{eq:lrv_estimation}
	\hat{\sigma}^2(t) = \sum_{i=1}^n w_{\tau_n}(t, i) \frac{(S_{i - m_n + 1 ,i} - S_{i + 1 ,i + m_n})^2}{2m_n},
\end{equation}
for $t\in [m_n / n, 1 - m_n/n]$, where $S_{i, j} = \sum_{k=i}^{j} X_{k, n}$, and extend the definition to $[0, 1]$, by defining $\hat{\sigma}^2(t) = \hat{\sigma}^2(m_n/n)$, for $t\in [0, m_n/n)$, and $\hat{\sigma}^2(t) = \hat{\sigma}^2(1 - m_n/n)$, for $t \in (1-m_n/n, 1]$.  If $m_n^{1/4} = o(\sqrt{n}\tau_n)$ and $m_n^{5/2} = o(n)$, then
\begin{equation*}
	\sup_{t\in[\gamma_n, 1-\gamma_n]} |\hat{\sigma}^2(t) - \sigma^2(t)| = o_\pr(1),
\end{equation*} 
where $\gamma_n = \tau_n + m_n/n$, by Theorem 6.4 of \cite{bucher2021}.

We are now able to define a test with asymptotic level $\alpha\in (0, 1)$. More specifically, we reject $H_0$ whenever 
\begin{equation} \label{eq:test_decision_linear}
	|S_n(\phi_1)| > \bigg(\int_0^1 \hat{\sigma}^2(t) \phi_1^2(t) \diff t\bigg)^{1/2} q_{1- \alpha/2},
\end{equation}
where $q_{1-\alpha/2}$ denotes the $(1-\alpha/2)$-quantile of the standard normal distribution. From the construction of the test and the expansion in \eqref{eq:lin_expansion}, we have the following proposition.

\begin{proposition} \label{prop:lin_test}
	Let $(h_n)_{n\in\N}$ be a bandwidth with $h_n = o(1)$ and $n^{1/4}h_n/\log^2 n \to \infty$. Under Assumptions \ref{assump:kern} and \ref{assump:error}, the test decision in \eqref{eq:test_decision_linear} defines a consistent test for the hypotheses in \eqref{eq:hypotheses_lin} with asymptotic level $\alpha$.
\end{proposition}

A crucial observation is that the hypotheses can be formulated via $\phi_1$ as 
\begin{equation} \label{eq:hypotheses_simple}
	H_0: \langle \phi_1, \mu - \bar{\mu} \rangle = 0 \quad \mathrm{vs.} \quad H_1: \langle \phi_1, \mu - \bar{\mu} \rangle \neq 0.
\end{equation}
where $\bar{\mu} = \int_0^1 \mu(x)\diff x$.
It turns out that the decision rule in \eqref{eq:test_decision_linear} does not depend on the choice of $\phi$, and $\phi_1$ can be replaced by any test function in $L^2([0, 1])$. In this case, the decision rule generalizes to
\begin{equation*} 
	|S_n(\phi)| > \bigg(\int_0^1 \hat{\sigma}^2(t) \Big(\phi(t) - \int_0^1 \phi(x)\diff x\Big)^2  \diff t\bigg)^{1/2} q_{1- \alpha/2}.
\end{equation*}

\subsection{Testing for a Polynomial Trend} \label{sec:poly_test}

We can extend the linear case to a polynomial mean function $\mu(t) = \alpha_0 + \sum_{i=1}^{p} \alpha_i t^i$, for some degree $p\in\N$. Then, the testing problem is equivalent to
\begin{equation}\label{eq:hypotheses_poly}
	H_0^{(\mathrm{poly})}: \alpha_i = 0~\mathrm{for~all}~i=1,\dots, p \quad \mathrm{vs.} \quad H_1^{(\mathrm{poly})}: \alpha_i \neq 0~\mathrm{for~some}~i.
\end{equation}
Similarly to the test function $\phi_1$, that can detect changes in a linear mean, we define the test functions $\phi_i(t) = t^i - \int_0^1 x^i \diff x$, for $i=1,\dots, p$. By Theorem \ref{thm:main}, the joint vector $(S_n(\phi_1), \dots, S_n(\phi_p))$ converges weakly to a centered normal distribution with covariance matrix $\Sigma = (\Sigma_{i,j})_{i, j=1}^p$, with $\Sigma_{i,j} = \int_0^1 \sigma^2(t) \phi_i(t) \phi_j(t) \diff t$ under the null hypothesis.

Based on the long-run variance estimator $\hat{\sigma}(t)$, we define a Wald-type test and reject $H_0^{(\mathrm{poly})}$, whenever
\begin{equation}  \label{eq:test_decision_poly}
	\Sb_n^T \hat{\Sigma}^{-1} \Sb_n > \chi_{p, 1-\alpha}^2,
\end{equation}
where $\Sb_n = (S_n(\phi_1), \dots, S_n(\phi_p))^T, \hat{\Sigma}_{i,j} = \int_0^1 \hat{\sigma}^2(t) \phi_i(t) \phi_j(t)$, and $\chi_{p, 1-\alpha}^2$ denotes the $(1-\alpha)$-quantile of the $\chi^2$-distribution with $p$ degrees of freedom.

\begin{proposition} \label{prop:poly_test}
	Let $(h_n)_{n\in\N}$ be a bandwidth with $h_n = o(1)$ and $n^{1/4}h_n/\log^2 n \to \infty$. Under Assumptions \ref{assump:kern} and \ref{assump:error}, the test decision in \eqref{eq:test_decision_poly} defines a consistent test for the hypotheses in \eqref{eq:hypotheses_poly} with asymptotic level $\alpha$.
\end{proposition}

As in Section \ref{sec:lin_test}, the test decision cannot only be used to test $H_0^{(\mathrm{poly})}$, but any composite hypothesis consisting of simple null hypotheses, as in \eqref{eq:hypotheses_simple}.

\subsection{Testing for Arbitrary Changes}

In the most general case, we want to test the hypotheses from \eqref{eq:hypotheses}, i.\,e., whether $\mu$ is constant against the alternative that it is not constant.
Using the approach from the previous section, we can only get consistency against alternatives in the span of finitely many test functions $\phi_1, \dots, \phi_p$. To get consistency against the general alternative $H_1$, we may consider a random test function $\Phi \in L^2([0, 1])$, and construct a test analogously to \eqref{eq:test_decision_linear}. 

More specifically, let $(e_k)_{k\in\N}$ be an orthonormal basis of $L^2([0, 1])$, such as the Fourier basis. Let $(\xi_k)_{k\in\N}$ be independent, standard normally distributed random variables, independent of the time series $X$, and let $(a_k)_{k\in\N}$ be a deterministic sequence with $\sum_{k\in\N} a_k^2 < \infty$. Define the random test function
\begin{equation*}
	\Phi(t) = \sum_{k=1}^{\infty} a_k \xi_k e_k(t),
\end{equation*}
which belongs a.s. to $L^2([0, 1])$ by Parseval's identity. Analogously to \eqref{eq:stat_lin}, define
\begin{equation*}
	S_n(\Phi) = \sqrt{n} \langle \Phi - \bar{\Phi}, \tilde{\mu}_{n}\rangle,
\end{equation*}
where $\bar{\Phi} = \int_0^1 \Phi(t) \diff t$, and subsequently, reject $H_0$ whenever
\begin{equation} \label{eq:test_decision_general} 
	|S_n(\Phi)| > \bigg(\int_0^1 \hat{\sigma}^2(t) \Phi^2(t) \diff t\bigg)^{1/2} q_{1- \alpha/2}.
\end{equation}

\begin{proposition} \label{prop:general_test}
	Let $(h_n)_{n\in\N}$ be a bandwidth with $h_n = o(1)$ and $n^{1/4}h_n/\log^2 n \to \infty$. Under Assumptions \ref{assump:kern} and \ref{assump:error}, the test decision in \eqref{eq:test_decision_general} defines a consistent test for the hypotheses in \eqref{eq:hypotheses} with asymptotic level $\alpha$.
\end{proposition}

The test defined by \eqref{eq:test_decision_general} is highly dependent on the realization of $\Phi$. For improved finite-sample robustness, we can use random test functions $\Phi_1, \dots, \Phi_p \in L^2([0, 1])$ with the decision rule in \eqref{eq:test_decision_poly}.

	% !TeX spellcheck = en_US
% ====================
\section{Empirical Results} \label{sec:empirical}
% ====================

We study the finite sample properties of the new methodology by means of a large scale simulation study and illustrate its application in a data example.\footnote{Python implementations of the methods and experiments are available on GitHub: \url{https://github.com/FlorianHeinrichs/Localized_FCLT}} For a comparative analysis, we use the bootstrap procedure proposed by \cite{zhou2013}. Moreover, we consider the classic CUSUM test as a benchmark. Under local stationarity, the partial sum process $n^{-1/2} \sum_{i=1}^{\lfloor nt\rfloor}( X_{i, n} - \ex[X_{i, n}])$ converges weakly to a centered Gaussian process $G$, defined by $G(t) = \int_0^t \sigma(x) \diff B_x$, for a standard Brownian motion $B$. To approximate the distribution of the CUSUM statistic, we simulate realizations of the estimated limiting process $\hat{G}(t) = \int_0^t \hat{\sigma}(x) \diff B_x$, where $\hat{\sigma}^2(t)$ is the local long-run variance estimator defined in \eqref{eq:lrv_estimation}.

For all experiments, we used decision rule \eqref{eq:test_decision_poly} with $p = 10$. Moreover, we constructed $p = 10$ random test functions using the first $100$ sine-cosine pairs of the Fourier basis with coefficients $a_k = \lfloor (k + 1) / 2\rfloor$, for $k=1, \dots, 200$. 

\subsection{Simulation Study} \label{sec:sim_study}

\begin{figure}
	\centering
	\includegraphics[width=\linewidth]{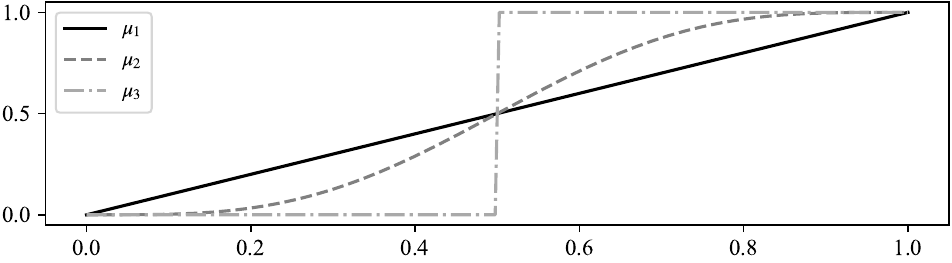}
	\caption{Various mean functions, used to generate time series under the alternative.}
	\label{fig:mu}
\end{figure}

Throughout, we consider the model
\begin{equation*}
	X_{i, n} = \mu(i/n) + \sigma(i/n) \eps_i,
\end{equation*}
for $i=1, \dots, n$, where $\mu$ denotes the mean function of interest, $\sigma^2(\cdot)$ a possibly time-varying variance and $(\eps_i)_{i\in\N}$ an error process.
The following choices of the mean function $\mu$ were considered
\[
\begin{array}{ll}
	\mu_0(x)  = 0, &  \mu_1(x)  = x, \\ 
	\mu_2(x)  = 35 x^4 - 84 x^5 + 70 x^6 - 20 x^7, & \mu_3(x)  = \id(x > \tfrac{1}{2}), 
\end{array}
\]
The function $\mu_0$ coincides with the null hypothesis of a constant mean, whereas the other functions induce a time-varying mean. The functions were selected to be linear, polynomial and non-polynomial, as displayed in Figure \ref{fig:mu}.
Similarly, for $\sigma^2$, we considered
\begin{equation*}
	\sigma_0(x)  = 1/2 \quad \text{and} \quad \sigma_1(x) =  1/2 - \cos(2 \pi x) / 4.
\end{equation*}
Finally, as error processes, a sequence of i.i.d. random variables and a locally stationary process were considered. More specifically, for independent sequences of i.i.d. random variables $(\eta_i^{(1)})_{i\in\Z}$ and $(\eta_i^{(2)})_{i\in\Z}$ with $\eta_i^{(1)} \sim \Nc(0, 1)$ and $\eta_i^{(2)} \sim \Uc([-\sqrt{3}, \sqrt{3}])$, we considered
\[ 
(\mathrm{iid})~ \eps_i = \eta_i^{(1)}, \qquad \text{and} \qquad (\mathrm{ls})~\eps_{i, n} = \sqrt{a(i/n)} \eps_i^{(1)} + \sqrt{1 - a(i/n)} \eps_i^{(2)},
\]
where $\eps_i^{(k)}$ satisfy the AR-equation $\eps_i^{(k)} = \tfrac{\sqrt{3}}{2}(\eta_i^{(k)} + \tfrac{1}{2} \eps_{i-1}^{(k)})$, for $k=1, 2$, and
\[
a(t) = \tfrac{1}{2}\big[1 - \cos\big(\tfrac{\pi}{2}[-\cos(\pi t) + 1]\big) \big].
\]

For all settings, we generated $1000$ time series and test $H_0$ with level $\alpha = 5\%$. Table \ref{tab:h0} contains  empirical rejection rates under the null hypothesis for different choices of $\eps, \sigma$ and $n$, while Table \ref{tab:h1} contains those values under the different alternatives. Additional results for MA- and AR-errors, as well as different choices of $\sigma$ can be found in Appendix \ref{app:empirical}.

First consider the null hypothesis $\mu = \mu_0$. For i.i.d. errors, the bootstrap- and CUSUM-tests have levels close to the nominal level of $\alpha = 0.05$. The test based on random projections, defined by \eqref{eq:test_decision_general}, seems rather conservative. Interestingly, the tests for linear and polynomial trends, defined by \eqref{eq:test_decision_linear} and \eqref{eq:test_decision_poly} seem conservative for time-varying $\sigma$, but have increasing empirical rejection rates exceeding the nominal level for constant $\sigma$. For locally stationary errors, all tests exceed the nominal level, where the bootstrap test is still closest to $\alpha = 0.05$. Interestingly, for constant $\sigma$, the proposed tests exceed the nominal level, whereas results are close to $5\%$ for time-varying $\sigma$. 

Under the alternative, the results are similar. The CUSUM test and \eqref{eq:test_decision_linear} seem powerful, and reject the null hypothesis in (almost) all cases even for small sample sizes. The bootstrap procedure seems to be slightly more conservative than the test based on random projections, \eqref{eq:test_decision_general}. The polynomial test yields results between \eqref{eq:test_decision_linear} and \eqref{eq:test_decision_general}.

\begin{table} \label{tab:h0}
	\caption{Empirical rejection rates for various choices of $\sigma$ and $\eps$ under the null hypothesis, for $\mu = \mu_0$}
	\begin{tabular}{llrrrrr}
		\toprule
		$n$ & $\sigma$ & Zhou2013 & CUSUM & \eqref{eq:test_decision_linear} & \eqref{eq:test_decision_poly} & \eqref{eq:test_decision_general} \\
		\midrule
		\multicolumn{7}{l}{\textit{Panel A: i.i.d. errors $\eps$}}\\		
		100 & $\sigma_0$ & 5.5 & 6.1 & 6.1 & 3.5 & 1.0 \\
		100 & $\sigma_1$ & 4.4 & 2.2 & 0.1 & 0.0 & 0.0 \\
		200 & $\sigma_0$ & 4.6 & 4.0 & 5.9 & 3.0 & 0.5 \\
		200 & $\sigma_1$ & 4.4 & 2.8 & 0.5 & 0.3 & 0.2 \\
		500 & $\sigma_0$ & 2.8 & 3.1 & 9.1 & 10.3 & 1.1 \\
		500 & $\sigma_1$ & 3.5 & 3.0 & 1.6 & 1.1 & 0.1 \\
		1000 & $\sigma_0$ & 3.8 & 5.7 & 13.9 & 28.6 & 3.4 \\
		1000 & $\sigma_1$ & 4.2 & 4.1 & 4.5 & 6.1 & 0.6 \\
		\midrule
		\multicolumn{7}{l}{\textit{Panel B: locally stationary errors $\eps$}}\\		
		100 & $\sigma_0$ & 8.4 & 28.7 & 16.2 & 24.3 & 11.0 \\
		100 & $\sigma_1$ & 6.4 & 13.8 & 3.1 & 6.0 & 2.8 \\
		200 & $\sigma_0$ & 7.6 & 25.8 & 17.1 & 21.6 & 9.9 \\
		200 & $\sigma_1$ & 4.9 & 15.6 & 4.0 & 6.1 & 4.4 \\
		500 & $\sigma_0$ & 5.7 & 20.9 & 14.2 & 28.3 & 11.1 \\
		500 & $\sigma_1$ & 4.4 & 14.6 & 6.2 & 5.9 & 5.1 \\
		1000 & $\sigma_0$ & 5.3 & 16.3 & 11.0 & 34.6 & 11.2 \\
		1000 & $\sigma_1$ & 4.8 & 13.8 & 7.4 & 11.6 & 5.5 \\
		\bottomrule
	\end{tabular}
\end{table}

\begin{table} \label{tab:h1}
	\caption{Empirical rejection rates for various choices of $\sigma$ and $\eps$ under various alternatives.}
	\begin{tabular}{lllrrrrr}
		\toprule
		$\eps$ & $n$ & $\sigma$ & Zhou2013 & CUSUM & \eqref{eq:test_decision_linear} & \eqref{eq:test_decision_poly} & \eqref{eq:test_decision_general} \\
		\midrule
		\multicolumn{8}{l}{\textit{Panel A: $\mu = \mu_1$}}\\		
		iid & 100 & $\sigma_0$ & 86.8 & 100.0 & 100.0 & 97.2 & 87.1 \\
		iid & 100 & $\sigma_1$ & 89.9 & 99.9 & 100.0 & 91.9 & 63.4 \\
		iid & 200 & $\sigma_0$ & 94.8 & 100.0 & 100.0 & 100.0 & 99.7 \\
		iid & 200 & $\sigma_1$ & 96.4 & 100.0 & 100.0 & 100.0 & 99.7 \\
		ls & 100 & $\sigma_0$ & 80.0 & 99.8 & 100.0 & 98.6 & 93.0 \\
		ls & 100 & $\sigma_1$ & 80.1 & 99.9 & 100.0 & 95.4 & 80.0 \\
		ls & 200 & $\sigma_0$ & 87.7 & 100.0 & 100.0 & 100.0 & 99.8 \\
		ls & 200 & $\sigma_1$ & 92.2 & 100.0 & 100.0 & 100.0 & 99.2 \\
		\midrule
		\multicolumn{8}{l}{\textit{Panel B: $\mu = \mu_2$}}\\		
		iid & 100 & $\sigma_0$ & 95.3 & 100.0 & 100.0 & 99.9 & 99.7 \\
		iid & 100 & $\sigma_1$ & 95.7 & 100.0 & 100.0 & 99.9 & 98.7 \\
		iid & 200 & $\sigma_0$ & 99.7 & 100.0 & 100.0 & 100.0 & 100.0 \\
		iid & 200 & $\sigma_1$ & 99.6 & 100.0 & 100.0 & 100.0 & 100.0 \\
		ls & 100 & $\sigma_0$ & 90.5 & 100.0 & 100.0 & 100.0 & 100.0 \\
		ls & 100 & $\sigma_1$ & 88.1 & 100.0 & 100.0 & 100.0 & 99.1 \\
		ls & 200 & $\sigma_0$ & 97.8 & 100.0 & 100.0 & 100.0 & 100.0 \\
		ls & 200 & $\sigma_1$ & 97.5 & 100.0 & 100.0 & 100.0 & 100.0 \\
		\midrule
		\multicolumn{8}{l}{\textit{Panel C: $\mu = \mu_3$}}\\		
		iid & 100 & $\sigma_0$ & 98.4 & 100.0 & 100.0 & 100.0 & 100.0 \\
		iid & 100 & $\sigma_1$ & 98.9 & 100.0 & 100.0 & 100.0 & 100.0 \\
		iid & 200 & $\sigma_0$ & 100.0 & 100.0 & 100.0 & 100.0 & 100.0 \\
		iid & 200 & $\sigma_1$ & 100.0 & 100.0 & 100.0 & 100.0 & 100.0 \\
		ls & 100 & $\sigma_0$ & 95.6 & 100.0 & 100.0 & 100.0 & 100.0 \\
		ls & 100 & $\sigma_1$ & 93.6 & 100.0 & 100.0 & 100.0 & 99.9 \\
		ls & 200 & $\sigma_0$ & 99.8 & 100.0 & 100.0 & 100.0 & 100.0 \\
		ls & 200 & $\sigma_1$ & 99.3 & 100.0 & 100.0 & 100.0 & 100.0 \\
		\bottomrule
	\end{tabular}
\end{table}

\subsection{Case Study} \label{sec:case_study}

In the following, we apply the proposed tests to detect changes in daily temperatures and brain activity.

\subsubsection{Daily Temperature}

We used temperature data from station measurements in three German cities (Aachen, Bochum and Göttingen) provided by the Deutscher Wetterdienst (DWD) Climate Data Center\footnote{\url{https://opendata.dwd.de/climate_environment/CDC/observations_germany/climate/10_minutes/air_temperature/}}. The dataset contains air temperature measured in 10-minute intervals. To handle the large number of missing values, we calculated the daily mean temperature. Days without any measurement were removed from the dataset. 
The dataset spans different time intervals in the three cities (Aachen: 1993 to 2011, Bochum: 2007 to 2024, Göttingen: 1993 to 2024). 

$p$-values of the tests are displayed in Table \ref{tab:temperature} for two time periods, the entire recording and a single year. For a single year, the $p$-values of all tests are below $1\%$, hence highly significant at a level of $5\%$. Similarly, the $p$-values of all tests for the entire history of Göttingen are below $1\%$. Only for the full records of Aachen and Bochum the test decisions differ. The bootstrap procedure by \cite{zhou2013} does not reject the null hypothesis of a change in both locations, whereas \eqref{eq:test_decision_linear} rejects only the null hypothesis in Bochum. The other tests reject all null hypotheses, indicating a statistically significant change in the temperature.

\begin{table} \label{tab:temperature}
	\caption{$p$-values (in \%, rounded to $1$ decimal) of tests for daily mean temperatures in Aachen (AC), Bochum (BO) and Göttingen (GÖ).}
	\begin{tabular}{l|rrr|rrr}
		\toprule
		& AC (full) & BO (full) & GÖ (full) & AC (1 year) & BO (1 year) & GÖ (1 year) \\
		\midrule
		\eqref{eq:test_decision_linear} & 7.1 & \textbf{0.0} & \textbf{0.0} & \textbf{0.0} & \textbf{0.7} & \textbf{0.1} \\
		\eqref{eq:test_decision_poly} & \textbf{0.0} & \textbf{0.0} & \textbf{0.0} & \textbf{0.0} & \textbf{0.0} & \textbf{0.0} \\
		\eqref{eq:test_decision_general} & \textbf{0.0} & \textbf{0.0} & \textbf{0.0} & \textbf{0.0} & \textbf{0.0} & \textbf{0.0} \\
		CUSUM & \textbf{0.1} & \textbf{0.0} & \textbf{0.0} & \textbf{0.0} & \textbf{0.0} & \textbf{0.0} \\
		Zhou2013 & 99.9 & 79.8 & \textbf{0.8} & \textbf{0.0} & \textbf{0.0} & \textbf{0.0} \\
		\bottomrule
	\end{tabular}
\end{table}

\subsubsection{EEG-Data Analysis}

We employed the ``Consumer-Grade EEG and Eye-Tracking Dataset``, containing EEG recordings from 116 sessions \cite{afonso2025}. The EEG was measured with 4 electrodes at a sampling rate of 256\,Hz, and we considered only the 102 sessions without technical difficulties. Each session includes four recordings, accumulating to 408 recordings.

As usual in neuroscience, we calculated a spectrogram of the raw time series with Hann windows of 2\,s and a step width of 0.5\,s, and applied the logarithm to the resulting spectra. The result is a functional time series, where at each step, we observe a spectrum based on 2\,s. To obtain univariate time series, we calculated the average power in the $\alpha$- and $\beta$-frequency bands, i.\,e.,
\begin{equation*}
	\alpha(f_i) = \frac{1}{12 - 8} \int_8^{12} f_i(x) \diff x, \quad \beta(f_i) = \frac{1}{30-13} \int_{13}^{30} f_i(x)  \diff x,
\end{equation*}
where $f_i$ denotes the spectrum of electrode $i \in\{1, \dots, 4\}$. We then tested each resulting time series for changes. Changes in these frequency bands are associated with a varying focus of participants.

The results are displayed in Table \ref{tab:eeg}. For each test, the number of recordings with detected changes is rather similar across the different electrodes, which indicates similar brain activity in the respective frequency bands. A notable exception is electrode TP9 in the $\beta$-frequency band, where substantially less changes were detected compared to other electrodes in the same band. Moreover, the tests \eqref{eq:test_decision_linear} and \eqref{eq:test_decision_general} yield similar results, whereas the CUSUM test and \eqref{eq:test_decision_poly} detect slightly more changes. As already indicated by the previous results, the bootstrap procedure appears to be more conservative, with substantially less detected changes; again TP9 in the $\beta$-frequency band is a notable exception. 

%\begin{table} \label{tab:eeg}
%	\caption{Ratio of recordings with a detected change (in \%). The column names (TP9, TP10, AF7 and AF8) refer to the electrodes' positions.}
%	\begin{tabular}{l|rrrrrrrr}
%		\toprule
%		& TP9 ($\alpha$) & TP9 ($\beta$) & AF7 ($\alpha$) & AF7 ($\beta$) & AF8 ($\alpha$) & AF8 ($\beta$) & TP10($\alpha$)  & TP10 ($\beta$) \\
%		\midrule
%		\eqref{eq:test_decision_linear} & 18.1 & 29.4 & 20.1 & 41.4 & 17.6 & 42.4 & 19.6 & 39.2 \\
%		\eqref{eq:test_decision_poly} & 22.1 & 40.2 & 27.0 & 62.0 & 28.7 & 63.5 & 25.7 & 50.2 \\
%		\eqref{eq:test_decision_general} & 14.7 & 29.9 & 15.0 & 41.7 & 18.1 & 41.9 & 17.4 & 38.5 \\
%		CUSUM & 20.1 & 33.6 & 25.0 & 49.8 & 21.3 & 48.3 & 22.3 & 47.3 \\
%		Zhou2013 & 8.8 & 10.3 & 7.6 & 15.4 & 7.1 & 17.4 & 8.3 & 13.7 \\
%		\bottomrule
%	\end{tabular}
%\end{table}

\begin{table} \label{tab:eeg}
	\caption{Ratio of recordings with a detected change (in \%). The column names (TP9, TP10, AF7 and AF8) refer to the electrodes' positions.}
	\begin{tabular}{l|rrrr|rrrr}
	\toprule
	& \multicolumn{4}{c|}{$\alpha$} & \multicolumn{4}{c}{$\beta$}\\
	& TP9 & AF7 & AF8 & TP10 & TP9 & AF7 & AF8 & TP10 \\
	\midrule
	\eqref{eq:test_decision_linear} 	& 18.1 & 20.1 & 17.6 & 19.6 & 10.3 & 41.4 & 42.4 & 39.2 \\
	\eqref{eq:test_decision_poly} 		& 22.1 & 27.0 & 28.7 & 25.7 & 33.6 & 62.0 & 63.5 & 50.2 \\
	\eqref{eq:test_decision_general} 	& 14.7 & 15.0 & 18.1 & 17.4 & 29.9 & 41.7 & 41.9 & 38.5 \\
	CUSUM 								& 20.1 & 25.0 & 21.3 & 22.3 & 40.2 & 49.8 & 48.3 & 47.3 \\
	Zhou2013 							&  8.8 &  7.6 &  7.1 &  8.3 & 29.4 & 15.4 & 17.4 & 13.7 \\
	\bottomrule
\end{tabular}
\end{table}

	% !TeX spellcheck = en_US

%% ====================
\section{Proofs} \label{sec:proofs}
%% ====================

%% ====================
\subsection{Proof of Theorem \ref{thm:main}} \label{sec:proof_main}
%% ====================

First, define $\eps_{i, n} = X_{i, n} - \ex[X_{i, n}]$. By Proposition 5 of \cite{zhou2013}, a sequence $(V_i)_{i\in\N}$ of independent, standard normally distributed random variables exists, on a possibly richer probability space, such that
\begin{equation*}
	\max_{j=1}^n \bigg| \sum_{i=1}^{j} \eps_{i, n} - \sum_{i=1}^j \sigma(\tfrac{i}{n}) V_i\bigg| = o_\pr(n^{1/4} \log^2 n).
\end{equation*}
By the same summation-by-parts argument as in the proof of Theorem A.1 of \cite{bucher2021}, for $t \in [h_n, 1-h_n]$, 
\begin{equation}\label{eq:gaussian_approx}
	\bigg| \sum_{i=1}^n \eps_{i, n} K_{h_n}(\tfrac{i}{n}-t) -  \sum_{i=1}^n \sigma(\tfrac{i}{n}) V_i K_{h_n}(\tfrac{i}{n}-t) \bigg| = o_\pr(n^{1/4} \log^2 n).
\end{equation}
Note that $\frac{1}{\sqrt{nh_n}} \sum_{i=1}^n \sigma(\tfrac{i}{n}) V_i K_{h_n}(\tfrac{i}{n}-t)$ is normally distributed with mean zero and variance
\begin{equation*}
	\frac{1}{nh_n} \sum_{i=1}^n \sigma^2(\tfrac{i}{n}) K_{h_n}(\tfrac{i}{n}-t) 
	\xrightarrow{n\to\infty} \sigma_K^2(t),
\end{equation*}
so that the first part of the theorem follows since $n^{-1/4}h_n^{-1/2}\log^2n \to 0$.
For the second part, let $\phi\in C_c^\infty(0, 1)$ be an arbitrary test function, and define $\phi(x) = 0$ for $x\notin(0, 1)$. By \eqref{eq:gaussian_approx},
\begin{align*}
	\langle \phi, G_n\rangle
	& = \frac{1}{\sqrt{n}h_n} \sum_{i=1}^n \eps_{i, n} \int_0^1 \phi(t) K_{h_n}(\tfrac{i}{n}-t) \diff t \\
	& = \frac{1}{\sqrt{n}h_n} \sum_{i=1}^n \sigma(\tfrac{i}{n}) V_i \int_0^1 \phi(t) K_{h_n}(\tfrac{i}{n}-t) \diff t + o_\pr(n^{-1/4} h_n^{-1} \log^2 n).
\end{align*}
The non-vanishing term on the right-hand side is normally distributed with mean $0$ and variance
\begin{equation} \label{eq:var}
	\Sigma_n^{(1)} := \frac{1}{nh_n^2} \sum_{i=1}^{n} \sigma^2(\tfrac{i}{n}) \Big(\int_0^1 \phi(t) K_{h_n}(\tfrac{i}{n}-t) \diff t\Big)^2.
\end{equation}
Since $K$ is supported on $[-1, 1]$, by substitution,
\begin{align*}
	\int_0^1 \phi(t) K_{h_n}(\tfrac{i}{n}-t) \diff t
	& = \int_{i/n-h_n}^{i/n+h_n} \phi(t) K_{h_n}(\tfrac{i}{n}-t) \diff t \\
	% = \int_{-h_n}^{h_n} \phi(s + \tfrac{i}{n}) K\big(\tfrac{s}{h_n}\big) \diff s
	&= h_n \int_{-1}^1 \phi(s h_n + \tfrac{i}{n}) K(s) \diff s
	= h_n \phi(\tfrac{i}{n}) + \Oc(h_n^2),
\end{align*}
uniformly for $t \in [h_n, 1- h_n]$. Hence, we can simplify the variance in \eqref{eq:var} to 
\begin{equation*}
	\Sigma_n^{(1)} 
	= \frac{1}{n} \sum_{i=\lceil nh_n\rceil}^{\lfloor n (1-h_n)\rfloor} \sigma^2(\tfrac{i}{n}) \phi^2(\tfrac{i}{n}) + \Oc(h_n)
	= \int_0^1 \sigma^2(x) \phi^2(x) \diff x + \Oc(h_n + n^{-1}).
\end{equation*}
In particular, $\langle \phi, G_n\rangle \convw \Nc(0, \int_0^1 \sigma^2(x) \phi^2(x) \diff x)$. 
Now, let $\phi_1, \dots, \phi_d \in C_c^\infty(0, 1)$ be test functions, and $\alpha_1, \dots, \alpha_d \in \R$. By the previous arguments, 
\begin{equation*}
	\sum_{i=1}^d \alpha_i \langle \phi_i, G_n \rangle 
	=  \Big\langle \sum_{i=1}^d \alpha_i \phi_i, G_n \Big\rangle 
	\convw \Nc(0, \Sigma_{i,j}),
\end{equation*}
where 
\begin{equation*}
	\Sigma_{i,j} = \int_0^1 \sigma^2(x) \Big(\sum_{i=1}^d \alpha_i \phi_i(x)\Big)^2 \diff x 
	= \sum_{i, j=1}^{d} \alpha_i \alpha_j \int_0^1 \sigma^2(x) \phi_i(x) \phi_j(x) \diff x.
\end{equation*}
Since the coefficients $\alpha_1, \dots, \alpha_d$ are arbitrary, by the Cramér-Wold device,
\begin{equation} \label{eq:fdd_conv}
	( \langle \phi_1, G_n \rangle, \dots, \langle \phi_d, G_n \rangle )
	\convw \Nc(0, \Sigma^{(d)})
\end{equation}
in $\R^d$, where $\Sigma^{(d)} = (\Sigma_{i,j})_{i, j=1}^d$. \qed

%By the Portmanteau theorem, $\ex[\exp(i t \langle G_n, \phi\angle)] \to \exp[-\tfrac{1}{2}\int_0^1 \sigma^2(x) \phi^2(x)\diff x]$, for any $t\in [0, 1]$ and $\phi \in C_c^\infty(0, 1)$.
%
%By identifying $G_n$ with the constant-path process $\Gb_n(t) = G_n$, $\phi(\Gb_n)$ is constant for every $\phi \in C_c^\infty(0, 1)$, hence it is tight in $D([0, 1], \R)$, since $\langle \phi, G_n\rangle$ converges weakly in $\R$. Mitoma's theorem \citep[see, e.\,g., Theorem 6.13 by][]{walsh1986} gives that $\Gb$ is tight in $D([0, 1], D'(0,1))$, hence $G = \Gb(0)$ is tight in $D'(0, 1)$.
%
% Together with \eqref{eq:fdd_conv}, Theorem 6.3 of \cite{fonseca2018} implies $\Gb_n \convw \Gb$ in $D_\infty(D'(0, 1))$, and therefore $G_n \convw G_n$ in $D'(0, 1)$.

%% ====================
\subsection{Proof of Theorem \ref{thm:isonormal}} \label{sec:proof_isonormal}
%% ====================

Let $f\in L^2([0, 1])$ be arbitrary. Since $C_c^\infty(0, 1)$ is dense in $L^2([0, 1])$, a sequence $(\phi_k)_{k\in\N}$ exists, such that $\|f - \phi_k\|_{L^2([0, 1])} \to 0$ \citep[see, e.\,g., Theorem 7.10 of][]{chacon2016}. 
By Theorem \ref{thm:main}, $\ex[(\langle \phi_k, G\rangle - \langle \phi_m, G\rangle)^2] = \|\sigma(\phi_k - \phi_m)\|_{L^2([0, 1])}^2 \to 0$, such that $(\langle \phi_k, G\rangle)_{k\in\N}$ is a Cauchy sequence in $L^2(\Omega)$, and converges by completeness. Define $W_\sigma(f) = L^2(\Omega)-\lim_{k\to\infty}\langle \phi_k, G\rangle$. The limit is well-defined, since all approximating sequences have the same asymptotic covariance. By construction, $\ex[W_\sigma(f)] = 0$, and $\ex[W_\sigma(f) W_\sigma(g)] = \langle \sigma f, \sigma g\rangle$, for all $f, g\in L^2([0, 1])$.

Now, observe that
\begin{align*}
	\ex[\|G_n\|_{L^2([0, 1])}^2]
	& = \frac{1}{nh_n^2} \ex\bigg[\int_0^1 \bigg(\sum_{i=1}^n(X_{i,n} - \ex[X_{i,n}]) K_{h_n}(\tfrac{i}{n}-t)\bigg)^2 \diff t\bigg] \\
	& = \frac{1}{nh_n^2} \sum_{i, j=1}^n \cov(X_{i,n}, X_{j,n}) \int_0^1 K_{h_n}(\tfrac{i}{n}-t) K_{h_n}(\tfrac{j}{n}-t) \diff t  \\
	% & = \frac{1}{nh_n} \sum_{i,j=1}^n \cov(X_{i,n}, X_{j,n}) \int_{-1}^1 K(-t) K(\tfrac{j-i}{nh_n}-t) \diff t \\
	& = \frac{1}{nh_n} \sum_{i=1}^n \sum_{j=(i-nh_n)\vee 1}^{(i+nh_n) \wedge n} \cov(X_{i,n}, X_{j,n}) \int_{-1}^1 K(-t) K(\tfrac{j-i}{nh_n}-t) \diff t,
\end{align*}
since $K$ has support $[-1, 1]$. Moreover, the integral on the right-hand side can be bounded by $\|K\|_{L^2([-1, 1])}^2$ by the Cauchy-Schwarz inequality. By Proposition 13 of \cite{heinrichs2025b}, $|\cov(X_{i,n}, X_{j,n})| \le \Theta_{|i-j|-1}$, for $\Theta_{i} = \sum_{k=i}^{\infty} \delta_2(k)$. In particular,
\begin{equation} \label{eq:moment_bound}
	\ex[\|G_n\|_{L^2([0, 1])}^2]
	\le \frac{1}{nh_n} \sum_{i=1}^n \sum_{j=(i-nh_n)\vee 1}^{(i+nh_n) \wedge n}  \Theta_{|i-j|-1} \|K\|_{L^2([-1, 1])}^2 \le C_K,
\end{equation}
for some constant $C_K > 0$, only depending on $K$ and the physical dependence measure $\delta_2(k)$.

Fix $f_1, \dots, f_d \in L^2([0, 1])$, and approximations $\{\phi_{i, m}\}_{i=1,\dots,d, m\in\N}\subset C_c^\infty(0, 1)$, such that $\| \phi_{i, m} - f_i\|_{L^2([0, 1])} \to 0$, as $m\to\infty$. Using the Cauchy-Schwarz inequality and \eqref{eq:moment_bound},
\begin{equation*}
	\ex[\langle f_i - \phi_{i, m}, G_n\rangle^2] \le C_K \| f_i - \phi_{i, m}\|_{L^2([0, 1])}^2,
\end{equation*}
which vanishes as $m \to\infty$. In particular, 
\begin{equation*}
	(\langle \phi_{1, m}, G_n \rangle, \dots, \langle \phi_{d, m}, G_n \rangle) \xrightarrow{L^2(\Omega)} (\langle f_1, G_n \rangle, \dots, \langle f_d, G_n \rangle).
\end{equation*}
Similarly, $\ex[(W_\sigma(f_i) - W_\sigma(\phi_{i, m}))^2] = \|\sigma(f_i - \phi_{i, m})\|_{L^2([0, 1])}^2 \to 0$, hence, 
\begin{equation*}
	(W_\sigma(\phi_{1, m}), \dots, W_\sigma(\phi_{d, m})) \xrightarrow{L^2(\Omega)} (W_\sigma(f_1), \dots, W_\sigma(f_d)).
\end{equation*}
Let $\theta\in\R^d$, and let $\phi_U(\theta) = \ex[\exp(i\theta^T U)]$ denote the characteristic function of some random vector $U\in\R^d$. For random vectors $U, V\in \R^d$, 
\begin{equation*}
	|\phi_U(\theta) - \phi_V(\theta)| \le \ex[|\exp(i\theta^T U) - \exp(i \theta^T V)|] \le \|\theta\|_2 \ex[\|U - V\|_2],
\end{equation*}
by the Cauchy-Schwarz inequality and $|\exp(ix)-\exp(iy)|\le |x-y|$, where $\|\cdot\|_2$ denotes the Euclidean norm in $\R^d$. Applying this inequality, yields
\begin{align*}
	& \sup_{n\in\N} |\phi_{(\langle f_1, G_n \rangle, \dots, \langle f_d, G_n \rangle)}(\theta) - \phi_{(\langle \phi_{1, m}, G_n \rangle, \dots, \langle \phi_{d, m}, G_n \rangle)}(\theta)| \to 0 \\
	& \sup_{n\in\N} |\phi_{(W_\sigma(f_1), \dots, W_\sigma(f_d))}(\theta) - \phi_{(W_\sigma(\phi_{1, m}), \dots, W_\sigma(\phi_{d, m}))}(\theta)| \to 0,
\end{align*}
as $m\to\infty$. By the Portmanteau theorem and Theorem \ref{thm:main}, for each fixed $m\in\N$, 
\begin{equation*}
	\phi_{(\langle \phi_{1, m}, G_n \rangle, \dots, \langle \phi_{d, m}, G_n \rangle)}(\theta) \xrightarrow{n\to\infty} \phi_{(W_\sigma(\phi_{1, m}), \dots, W_\sigma(\phi_{d, m}))}(\theta).
\end{equation*}
Combining the last three convergences, we finally have
\begin{equation*}
	\phi_{(\langle f_1, G_n \rangle, \dots, \langle f_d, G_n \rangle)}(\theta) \xrightarrow{n\to\infty} \phi_{(W_\sigma(f_1), \dots, W_\sigma(f_d))}(\theta),
\end{equation*}
so that the theorem's statement follows from the Portmanteau theorem.

% For each $n\in\N$, define the linear map $T_n: C_c^\infty(0, 1)\to L^2(\Omega)$ by $T_n(\phi) = \langle \phi, G_n\rangle$. By \eqref{eq:moment_bound}, each $T_n$ is bounded with operator norm bounded by $C_K$. Since $C_c^\infty(0, 1)$ is dense in $L^2([0, 1])$, $T_n$ extends uniquely by continuity to a bounded linear operator on $L^2([0, 1])$ by the Hahn-Banach theorem \citep[see, e.\,g., Corollary 7.3.3 in][]{narici2011}, such that $T_n(\phi) =  \langle \phi, G_n\rangle$ for $\phi \in L^2([0, 1])$.

%% ====================
\subsection{Proof of Theorem \ref{thm:hoelder}} \label{sec:proof_hoelder}
%% ====================

% Since $\Gc$ is totally bounded, it is separable, and the Gaussian limit process admits a version with uniformly $L^2$-continuous sample paths because its canonical semimetric is dominated by $C \|\cdot\|_{L^2([0, 1])}$.
Let $W_\sigma$ denote the Gaussian process from Theorem \ref{thm:isonormal}, and let $H:= W_\sigma\vert_\Gc$ denote the restriction of $W_\sigma$ to $\Gc$. By Theorem 1.5.4 of \cite{vandervaart2023}, $H_n$ converges weakly to $H$, if the finite-dimensional distributions converge and $H_n$ is asymptotically tight. By Theorem \ref{thm:isonormal}, the finite-dimensional distributions converge, i.\,e.,
\begin{equation*}
	(H_n(f_1), \dots, H_n(f_d)) \convw (H(f_1), \dots, H(f_d)),
\end{equation*}
for any $f_1,\dots, f_d\in \Gc$. By Theorem 1.5.7 of \cite{vandervaart2023}, $H_n$ is asymptotically tight, iff $H_n(f)$ is asymptotically tight for each $f\in\Gc$, $\Gc$ is totally bounded w.r.t. some seminorm $\rho$, and $H_n$ is asymptotically uniformly $\rho$-equicontinuous in probability. Asymptotic tightness of $H_n(f)$ in $\R$ follows immediately from its weak convergence. $(\Gc, \|\cdot\|_{L^2([0, 1])})$ is totally bounded by assumption. Finally, recall from \eqref{eq:moment_bound}, $\ex[\|G_n\|_{L^2([0, 1])}^2] \le C_K$, for some constant $C_K > 0$. Hence, for any $f, g \in \Gc$,
\begin{equation*}
	E[(H_n(f) - H_n(g))^2] % = \ex[\langle f-g, G_n \rangle^2] 
	\le \ex[\|G_n\|_{L^2([0, 1])}^2] \|f-g\|_{L^2([0, 1])}^2 \le C_K \|f-g\|_{L^2([0, 1])}^2,
\end{equation*}
by the Cauchy-Schwarz inequality. In particular, 
\begin{equation*}
	\pr\big(\sup_{\|f-g\|_{L^2}< \delta} |H_n(f) - H_n(g)| > \eps \big)
	\le \pr\big(\delta \|G_n\|_{L^2([0, 1])}^2 > \eps \big)
	\le \frac{C_K \delta^2}{\eps^2}
\end{equation*}
by Markov's inequality, and asymptotic uniform $L^2([0, 1])$-equicontinuity in probability follows.

\subsection{Proofs of the Results from Section \ref{sec:stat_app}}

% \subsubsection{Proof of Proposition \ref{prop:lin_test}}

\subsubsection{Proof of Proposition \ref{prop:poly_test}}

The covariance matrix $\Sigma$ is symmetric by definition. Moreover, the test functions $\phi_i(t)$ are linearly independent and the long-run variance is bounded away from zero by Assumption \ref{assump:error} (4), hence $\Sigma$ is positive definite and invertible.

Under the null hypothesis, $\Sb_n^T \hat{\Sigma}^{-1} \Sb_n$ converges weakly to $Y^T \Sigma^{-1} Y$, for $Y \sim \Nc(0, \Sigma)$, by Theorem \ref{thm:isonormal} and the continuous mapping theorem. By standard arguments based on the Cholesky decomposition of $\Sigma$, $\Sigma^{-1/2} Y \sim \Nc(0, I_p)$, where $I_p$ denotes the $p\times p$ identity matrix. Therefore, $Y^T \Sigma^{-1} Y \sim \chi_p^2$ and the test defined by \eqref{eq:test_decision_poly} has asymptotic level $\alpha$.

Conversely, under the alternative at least some index $\alpha_i$ is nonzero. Analogously to \eqref{eq:lin_expansion2}, 
%\begin{align*}
%	\sqrt{n} \int_0^1 \phi_i(t) \Big(\mu(t) - \int_0^1 \mu(x)\diff x\Big) \diff t
%	& = \sqrt{n} \sum_{k=1}^p \alpha_k  \int_0^1 \phi_i(t) \phi_k(t) \diff t \\
%	% & = \sqrt{n} \sum_{k=1}^p \alpha_k \bigg(\frac{1}{i+k+1} - \frac{1}{(i+1)(k+1)}\bigg)
%	& = \sqrt{n} \sum_{k=1}^p \alpha_k \frac{ik}{(i+1)(k+1)(i+k+1)},
%\end{align*}
\begin{equation*}
	\sqrt{n} \Big\langle \phi_i, \mu - \int_0^1 \mu(x)\diff x \Big\rangle
	= \sqrt{n} \sum_{k=1}^p \alpha_k  \langle \phi_i, \phi_k\rangle
\end{equation*}
Let $A = (A_{j,k})_{j,k=1}^p$ with $A_{j,k} = \langle \phi_j, \phi_k \rangle$. Then the drift of $\mu$ in direction $\phi_i$ is the $i$-th coordinate of $\sqrt{n}A \alpha$, where $\alpha = (\alpha_1, \dots, \alpha_p)^T$. Analogously to \eqref{eq:lin_expansion}, we may expand $\Sb_n = \sqrt{n} A \alpha + \Sb_n^c$, where $\Sb_n^c = \Sb_n - \sqrt{n} A \alpha$. By the same arguments as under the null hypothesis, $\Sb_n^c$ converges weakly to $\Nc(0, \Sigma)$. Hence, whenever $A \alpha \neq 0$, the Wald statistic $\Sb_n^T \hat{\Sigma}^{-1} \Sb_n$ is dominated by $n \alpha^T A^T \hat{\Sigma}^{-1} A \alpha$. 

Now observe that $A$ is the Gram matrix of the functions $\phi_1, \dots, \phi_p$ in $L^2([0, 1])$. For $\alpha\in\R^p$, define $h_\alpha(t) = \sum_{i=1}^p \alpha_i \phi_i(t)$. By linearity of the scalar product,
\begin{equation*}
	\alpha^T A \alpha
	= \sum_{i, j=1}^p \alpha_i \alpha_j A_{i, j} = \langle h_\alpha, h_\alpha\rangle = \|h_\alpha\|_{L^2([0, 1])}^2 \ge 0.
\end{equation*}

Hence, $\alpha^T A \alpha = 0$ if and only if $h_\alpha(t) = 0$ almost everywhere. Since the functions $\phi_1, \dots, \phi_p$ are linearly independent, $\alpha \neq 0$ implies $h_\alpha(t) = \sum_{i=1}^p \alpha_i \phi_i(t) \not\equiv 0$. In particular, $A$ is symmetric positive definite and $\alpha \neq 0$ implies $A \alpha \neq 0$. Therefore, the Wald statistic diverges to $\infty$ under the alternative.

\subsubsection{Proof of Proposition \ref{prop:general_test}}

First, $\mu(t)$ is constant iff $\mu(t) = \int_0^1 \mu(x) \diff x$. Let $g(t) = \mu(t) - \int_0^1 \mu(x) \diff x$, so that the null hypothesis translates to $H_0: g(t) = 0$, for all $t\in[0, 1]$. By the expansion from \eqref{eq:lin_expansion}, 
\begin{equation*}
	S_n(\Phi) = \sqrt{n} \langle \Phi, \tilde{\mu}_{n} - \mu \rangle + \sqrt{n} \langle \Phi, g \rangle.
\end{equation*}
Since $\Phi$ is independent of $X$ and $\Phi \in L^2([0,1])$ almost surely, we argue conditionally on $\Phi$. Theorem \ref{thm:isonormal} yields
\begin{equation*}
	\sqrt{n} \langle \phi, \tilde{\mu}_n - \mu\rangle \convw \Nc \Big(0, \int_0^1 \sigma^2(t)\phi^2(t) \diff t \Big),
\end{equation*}
for $\pr_\Phi$-almost every realization $\Phi=\phi$. Equivalently,
\begin{equation*}
	\sqrt{n} \langle \Phi, \tilde{\mu}_n - \mu\rangle \vert \Phi \convw \Nc \Big(0, \int_0^1 \sigma^2(t)\Phi^2(t) \diff t \Big) \quad \text{a.s.}
\end{equation*}
Under the null hypothesis, $\langle \Phi, g \rangle = 0$, since $g \equiv 0$, and the test  has conditional asymptotic level $\alpha$ given $\Phi$, and thus also unconditional asymptotic level $\alpha$.

For any $g\in L^2([0, 1])$, $\langle \Phi, g\rangle = \int_0^1 \Phi(t) g(t) \diff t$ is a centered Gaussian random variable with variance
\begin{equation*}
	\var(\langle \Phi, g\rangle) = \sum_{k\in\N} a_k^2 \langle e_k, g\rangle^2.
\end{equation*}
Under the alternative, $g \not \equiv 0$, hence, $\langle e_k, g\rangle \neq 0$, for at least one $k\in\N$, such that the variance is strictly positive. Hence, $\pr(\langle \Phi, g \rangle = 0) = 0$. Therefore, for almost every realization of $\Phi$, the deterministic term $\sqrt n \langle \Phi,g\rangle$ diverges in absolute value, while the centered stochastic term remains conditionally $O_\pr(1)$ by the previous conditional weak convergence. This implies conditional consistency given $\Phi$ almost surely, and therefore unconditional consistency.

	\backmatter
	
%	\bmhead{Supplementary information}
%	
%	If your article has accompanying supplementary file/s please state so here. 
%	
%	Authors reporting data from electrophoretic gels and blots should supply the full unprocessed scans for key as part of their Supplementary information. This may be requested by the editorial team/s if it is missing.
%	
%	Please refer to Journal-level guidance for any specific requirements.
%	
%	\bmhead{Acknowledgements}
%	
%	Acknowledgements are not compulsory. Where included they should be brief. Grant or contribution numbers may be acknowledged.
%	
%	Please refer to Journal-level guidance for any specific requirements.
	
	\section*{Competing Interests}

	The authors have no competing interests to declare that are relevant to the content of this article.

%	\section*{Data Availability}

	%\bibliographystyle{apalike}
	\bibliography{bibliography}
	
	\begin{appendices}
		
		% !TeX spellcheck = en_US
% ====================
 \section{Additional Empirical Results} \label{app:empirical}
% ====================

This section contains additional results for the experiments from Section \ref{sec:sim_study} under the model $X_{i, n} = \mu(i/n) + \sigma(i/n) \eps_i$. In addition to the choices $\sigma_0$ and $\sigma_1$, we considered
\begin{equation*}
	\sigma_2(x)  = \tfrac{1}{4} + \tfrac{x}{2} \quad \text{and} \quad \sigma_3(x) =   \tfrac{1}{4} +  \tfrac{1}{2} \id(x >  \tfrac{1}{2}).
\end{equation*} 
Moreover, we considered moving average and autoregressive processes for $\eps$, i.\,e., 
\begin{equation*}
(\mathrm{ma})~ \eps_i = \tfrac{2}{\sqrt{5}}(\eta_i + \tfrac{1}{2} \eta_{i-1}) \qquad \text{and} \qquad (\mathrm{ar})~ \eps_i = \tfrac{\sqrt{3}}{2}(\eta_i + \tfrac{1}{2} \eps_{i-1}),
\end{equation*}
for $(\eta_i)_{i\in\Z}$ with $\eta_i\sim\Nc(0, 1)$ i.i.d.

\begin{table}
	\caption{Empirical rejection rates for various choices of $\sigma$ and $\eps$, under the null hypothesis $\mu = \mu_0$}
	\begin{tabular}{ll|rrrrr|rrrrr}
		\toprule
		$n$ & $\sigma$ & Zhou2013 & CUSUM & \eqref{eq:test_decision_linear} & \eqref{eq:test_decision_poly} & \eqref{eq:test_decision_general} & Zhou2013 & CUSUM & \eqref{eq:test_decision_linear} & \eqref{eq:test_decision_poly} & \eqref{eq:test_decision_general} \\
		\midrule
		&& \multicolumn{5}{l}{\textit{Panel A: i.i.d. errors $\eps$}} & \multicolumn{5}{l}{\textit{Panel B: MA errors $\eps$}}\\
		100 & $\sigma_0$ & 5.5 & 6.1 & 6.1 & 3.5 & 1.0  & 6.3 & 11.5 & 9.6 & 13.7 & 5.6 \\
		100 & $\sigma_1$ & 4.4 & 2.2 & 0.1 & 0.0 & 0.0  & 4.1 & 6.5 & 0.5 & 2.2 & 1.8 \\
		100 & $\sigma_2$ & 7.9 & 7.3 & 6.9 & 2.6 & 0.6  & 9.1 & 16.6 & 11.5 & 15.1 & 7.0 \\
		100 & $\sigma_3$ & 11.6 & 7.1 & 7.0 & 2.3 & 0.9  & 10.0 & 15.4 & 10.4 & 10.5 & 6.2 \\
		200 & $\sigma_0$ & 4.6 & 4.0 & 5.9 & 3.0 & 0.5  & 5.3 & 10.2 & 7.9 & 11.6 & 7.9 \\
		200 & $\sigma_1$ & 4.4 & 2.8 & 0.5 & 0.3 & 0.2  & 3.2 & 6.7 & 1.5 & 1.8 & 1.8 \\
		200 & $\sigma_2$ & 8.3 & 8.8 & 7.8 & 4.1 & 1.1  & 6.4 & 13.6 & 10.9 & 11.4 & 7.4 \\
		200 & $\sigma_3$ & 9.0 & 5.9 & 6.5 & 4.2 & 0.9  & 9.4 & 12.3 & 10.4 & 10.3 & 6.5 \\
		500 & $\sigma_0$ & 2.8 & 3.1 & 9.1 & 10.3 & 1.1  & 3.9 & 8.5 & 8.4 & 12.1 & 6.9 \\
		500 & $\sigma_1$ & 3.5 & 3.0 & 1.6 & 1.1 & 0.1  & 3.3 & 6.3 & 3.0 & 2.3 & 3.9 \\
		500 & $\sigma_2$ & 6.5 & 6.1 & 11.4 & 10.6 & 1.5  & 4.9 & 11.3 & 9.8 & 11.0 & 8.9 \\
		500 & $\sigma_3$ & 8.5 & 5.4 & 8.9 & 12.0 & 0.8  & 5.9 & 9.8 & 9.6 & 12.1 & 6.3 \\
		1000 & $\sigma_0$ & 3.8 & 5.7 & 13.9 & 28.6 & 3.4  & 3.9 & 7.9 & 7.8 & 13.2 & 6.0 \\
		1000 & $\sigma_1$ & 4.2 & 4.1 & 4.5 & 6.1 & 0.6  & 3.3 & 5.2 & 2.9 & 2.6 & 3.2 \\
		1000 & $\sigma_2$ & 4.4 & 5.8 & 13.7 & 22.3 & 3.9  & 5.1 & 9.4 & 8.4 & 11.8 & 7.7 \\
		1000 & $\sigma_3$ & 6.5 & 5.8 & 12.5 & 28.3 & 2.7  & 6.4 & 9.2 & 8.0 & 13.8 & 6.6 \\
		\midrule
		&& \multicolumn{5}{l}{\textit{Panel C: AR errors $\eps$}} & \multicolumn{5}{l}{\textit{Panel D: locally stationary errors $\eps$}}\\
		100 & $\sigma_0$ & 5.0 & 33.0 & 21.7 & 48.8 & 32.8  & 8.4 & 28.7 & 16.2 & 24.3 & 11.0 \\
		100 & $\sigma_1$ & 2.7 & 19.9 & 4.2 & 19.7 & 16.6  & 6.4 & 13.8 & 3.1 & 6.0 & 2.8 \\
		100 & $\sigma_2$ & 6.9 & 31.6 & 21.3 & 45.7 & 30.3  & 14.6 & 39.0 & 23.5 & 31.0 & 17.1 \\
		100 & $\sigma_3$ & 9.2 & 32.2 & 22.0 & 39.3 & 25.4  & 15.0 & 33.8 & 22.6 & 24.2 & 15.4 \\
		200 & $\sigma_0$ & 4.5 & 24.4 & 14.3 & 46.7 & 35.4  & 7.6 & 25.8 & 17.1 & 21.6 & 9.9 \\
		200 & $\sigma_1$ & 3.2 & 18.1 & 5.2 & 15.6 & 17.5  & 4.9 & 15.6 & 4.0 & 6.1 & 4.4 \\
		200 & $\sigma_2$ & 4.7 & 27.9 & 17.2 & 45.7 & 34.4  & 9.6 & 32.8 & 20.5 & 28.3 & 19.9 \\
		200 & $\sigma_3$ & 8.3 & 26.5 & 18.0 & 43.7 & 29.8  & 10.0 & 28.6 & 18.2 & 23.5 & 18.0 \\
		500 & $\sigma_0$ & 4.0 & 16.8 & 12.5 & 35.0 & 28.9  & 5.7 & 20.9 & 14.2 & 28.3 & 11.1 \\
		500 & $\sigma_1$ & 2.6 & 12.9 & 5.2 & 11.0 & 18.0  & 4.4 & 14.6 & 6.2 & 5.9 & 5.1 \\
		500 & $\sigma_2$ & 5.2 & 21.6 & 15.2 & 32.9 & 32.4  & 7.8 & 24.8 & 16.9 & 23.1 & 20.1 \\
		500 & $\sigma_3$ & 6.0 & 18.4 & 12.8 & 35.2 & 29.1  & 8.3 & 20.3 & 14.6 & 21.4 & 20.4 \\
		1000 & $\sigma_0$ & 4.4 & 14.7 & 11.1 & 24.3 & 24.7  & 5.3 & 16.3 & 11.0 & 34.6 & 11.2 \\
		1000 & $\sigma_1$ & 3.2 & 11.6 & 4.6 & 9.0 & 17.2  & 4.8 & 13.8 & 7.4 & 11.6 & 5.5 \\
		1000 & $\sigma_2$ & 4.9 & 16.6 & 12.8 & 21.6 & 23.4  & 7.5 & 19.1 & 15.3 & 20.1 & 20.2 \\
		1000 & $\sigma_3$ & 5.3 & 15.1 & 11.7 & 24.3 & 21.6  & 5.6 & 15.3 & 10.2 & 17.1 & 18.8 \\
		\bottomrule
	\end{tabular}
\end{table}

\begin{table}
	\caption{Empirical rejection rates for various choices of $\sigma$ and $\eps$, under the alternative $\mu = \mu_1$ }
	\begin{tabular}{ll|rrrrr|rrrrr}
		\toprule
		$n$ & $\sigma$ & Zhou2013 & CUSUM & \eqref{eq:test_decision_linear} & \eqref{eq:test_decision_poly} & \eqref{eq:test_decision_general} & Zhou2013 & CUSUM & \eqref{eq:test_decision_linear} & \eqref{eq:test_decision_poly} & \eqref{eq:test_decision_general}\\
		\midrule
		&&\multicolumn{5}{l}{\textit{Panel A: i.i.d. errors $\eps$}} & \multicolumn{5}{|l}{\textit{Panel B: MA errors $\eps$}}\\
		100 & $\sigma_0$ & 86.8 & 100.0 & 100.0 & 97.2 & 87.1  & 74.4 & 98.9 & 99.8 & 91.7 & 75.8 \\
		100 & $\sigma_1$ & 89.9 & 99.9 & 100.0 & 91.9 & 63.4  & 77.9 & 98.3 & 99.9 & 76.7 & 54.9 \\
		100 & $\sigma_2$ & 82.7 & 100.0 & 100.0 & 97.3 & 86.4  & 71.6 & 97.8 & 99.5 & 89.3 & 74.5 \\
		100 & $\sigma_3$ & 85.4 & 98.7 & 99.8 & 97.3 & 79.0  & 70.6 & 94.8 & 99.0 & 89.6 & 68.3 \\
		200 & $\sigma_0$ & 94.8 & 100.0 & 100.0 & 100.0 & 99.7  & 87.6 & 100.0 & 100.0 & 99.7 & 97.3 \\
		200 & $\sigma_1$ & 96.4 & 100.0 & 100.0 & 100.0 & 99.7  & 91.1 & 100.0 & 100.0 & 100.0 & 96.8 \\
		200 & $\sigma_2$ & 94.4 & 100.0 & 100.0 & 100.0 & 99.6  & 85.8 & 100.0 & 100.0 & 100.0 & 98.0 \\
		200 & $\sigma_3$ & 94.1 & 100.0 & 100.0 & 100.0 & 99.6  & 85.2 & 99.9 & 100.0 & 100.0 & 98.5 \\
		500 & $\sigma_0$ & 100.0 & 100.0 & 100.0 & 100.0 & 100.0  & 97.3 & 100.0 & 100.0 & 100.0 & 100.0 \\
		500 & $\sigma_1$ & 99.9 & 100.0 & 100.0 & 100.0 & 100.0  & 99.4 & 100.0 & 100.0 & 100.0 & 100.0 \\
		500 & $\sigma_2$ & 99.7 & 100.0 & 100.0 & 100.0 & 100.0  & 98.4 & 100.0 & 100.0 & 100.0 & 100.0 \\
		500 & $\sigma_3$ & 99.4 & 100.0 & 100.0 & 100.0 & 100.0  & 97.0 & 100.0 & 100.0 & 100.0 & 100.0 \\
		1000 & $\sigma_0$ & 100.0 & 100.0 & 100.0 & 100.0 & 100.0  & 100.0 & 100.0 & 100.0 & 100.0 & 100.0 \\
		1000 & $\sigma_1$ & 100.0 & 100.0 & 100.0 & 100.0 & 100.0  & 100.0 & 100.0 & 100.0 & 100.0 & 100.0 \\
		1000 & $\sigma_2$ & 100.0 & 100.0 & 100.0 & 100.0 & 100.0  & 99.8 & 100.0 & 100.0 & 100.0 & 100.0 \\
		1000 & $\sigma_3$ & 100.0 & 100.0 & 100.0 & 100.0 & 100.0  & 100.0 & 100.0 & 100.0 & 100.0 & 100.0 \\
		\midrule
		&&\multicolumn{5}{l}{\textit{Panel C: AR errors $\eps$}} & \multicolumn{5}{|l}{\textit{Panel D: locally stationary errors $\eps$}}\\
		100 & $\sigma_0$ & 57.9 & 98.3 & 98.9 & 95.6 & 88.1  & 80.0 & 99.8 & 100.0 & 98.6 & 93.0 \\
		100 & $\sigma_1$ & 62.5 & 97.6 & 99.3 & 88.2 & 74.1  & 80.1 & 99.9 & 100.0 & 95.4 & 80.0 \\
		100 & $\sigma_2$ & 55.4 & 96.8 & 97.1 & 95.2 & 85.7  & 69.0 & 98.8 & 99.5 & 98.2 & 91.4 \\
		100 & $\sigma_3$ & 56.4 & 93.4 & 95.5 & 91.8 & 81.0  & 64.0 & 95.1 & 98.0 & 96.3 & 84.5 \\
		200 & $\sigma_0$ & 75.7 & 99.9 & 99.9 & 99.6 & 97.6  & 87.7 & 100.0 & 100.0 & 100.0 & 99.8 \\
		200 & $\sigma_1$ & 80.4 & 99.9 & 100.0 & 99.6 & 96.4  & 92.2 & 100.0 & 100.0 & 100.0 & 99.2 \\
		200 & $\sigma_2$ & 73.5 & 100.0 & 99.9 & 99.3 & 96.7  & 80.8 & 100.0 & 100.0 & 100.0 & 99.6 \\
		200 & $\sigma_3$ & 68.6 & 99.2 & 99.8 & 99.7 & 97.2  & 73.7 & 99.4 & 99.9 & 100.0 & 98.1 \\
		500 & $\sigma_0$ & 94.9 & 100.0 & 100.0 & 100.0 & 100.0  & 98.8 & 100.0 & 100.0 & 100.0 & 100.0 \\
		500 & $\sigma_1$ & 97.4 & 100.0 & 100.0 & 100.0 & 100.0  & 100.0 & 100.0 & 100.0 & 100.0 & 100.0 \\
		500 & $\sigma_2$ & 94.7 & 100.0 & 100.0 & 100.0 & 100.0  & 94.5 & 100.0 & 100.0 & 100.0 & 100.0 \\
		500 & $\sigma_3$ & 92.1 & 100.0 & 100.0 & 100.0 & 100.0  & 93.5 & 100.0 & 100.0 & 100.0 & 100.0 \\
		1000 & $\sigma_0$ & 99.7 & 100.0 & 100.0 & 100.0 & 100.0  & 100.0 & 100.0 & 100.0 & 100.0 & 100.0 \\
		1000 & $\sigma_1$ & 99.9 & 100.0 & 100.0 & 100.0 & 100.0  & 100.0 & 100.0 & 100.0 & 100.0 & 100.0 \\
		1000 & $\sigma_2$ & 99.5 & 100.0 & 100.0 & 100.0 & 100.0  & 99.4 & 100.0 & 100.0 & 100.0 & 100.0 \\
		1000 & $\sigma_3$ & 99.3 & 100.0 & 100.0 & 100.0 & 100.0  & 99.3 & 100.0 & 100.0 & 100.0 & 100.0 \\
		\bottomrule
	\end{tabular}
\end{table}

\begin{table}
	\caption{Empirical rejection rates for various choices of $\sigma$ and $\eps$, under the alternative $\mu = \mu_2$ }
	\begin{tabular}{ll|rrrrr|rrrrr}
		\toprule
		$n$ & $\sigma$ & Zhou2013 & CUSUM & \eqref{eq:test_decision_linear} & \eqref{eq:test_decision_poly} & \eqref{eq:test_decision_general} & Zhou2013 & CUSUM & \eqref{eq:test_decision_linear} & \eqref{eq:test_decision_poly} & \eqref{eq:test_decision_general} \\
		\midrule
		&& \multicolumn{5}{l}{\textit{Panel A: i.i.d. errors $\eps$}} & \multicolumn{5}{|l}{\textit{Panel B: MA errors $\eps$}}\\
		100 & $\sigma_0$ & 95.3 & 100.0 & 100.0 & 99.9 & 99.7  & 85.3 & 100.0 & 100.0 & 99.3 & 97.6 \\
		100 & $\sigma_1$ & 95.7 & 100.0 & 100.0 & 99.9 & 98.7  & 84.8 & 100.0 & 100.0 & 97.8 & 92.6 \\
		100 & $\sigma_2$ & 95.1 & 100.0 & 100.0 & 100.0 & 99.7  & 82.7 & 100.0 & 100.0 & 99.3 & 97.9 \\
		100 & $\sigma_3$ & 94.5 & 100.0 & 100.0 & 99.9 & 99.5  & 82.5 & 99.9 & 99.9 & 98.9 & 95.6 \\
		200 & $\sigma_0$ & 99.7 & 100.0 & 100.0 & 100.0 & 100.0  & 93.2 & 100.0 & 100.0 & 100.0 & 100.0 \\
		200 & $\sigma_1$ & 99.6 & 100.0 & 100.0 & 100.0 & 100.0  & 97.1 & 100.0 & 100.0 & 100.0 & 100.0 \\
		200 & $\sigma_2$ & 99.6 & 100.0 & 100.0 & 100.0 & 100.0  & 94.1 & 100.0 & 100.0 & 100.0 & 100.0 \\
		200 & $\sigma_3$ & 99.4 & 100.0 & 100.0 & 100.0 & 100.0  & 92.0 & 100.0 & 100.0 & 100.0 & 99.9 \\
		500 & $\sigma_0$ & 100.0 & 100.0 & 100.0 & 100.0 & 100.0  & 99.7 & 100.0 & 100.0 & 100.0 & 100.0 \\
		500 & $\sigma_1$ & 100.0 & 100.0 & 100.0 & 100.0 & 100.0  & 99.9 & 100.0 & 100.0 & 100.0 & 100.0 \\
		500 & $\sigma_2$ & 100.0 & 100.0 & 100.0 & 100.0 & 100.0  & 99.9 & 100.0 & 100.0 & 100.0 & 100.0 \\
		500 & $\sigma_3$ & 100.0 & 100.0 & 100.0 & 100.0 & 100.0  & 99.6 & 100.0 & 100.0 & 100.0 & 100.0 \\
		1000 & $\sigma_0$ & 100.0 & 100.0 & 100.0 & 100.0 & 100.0  & 100.0 & 100.0 & 100.0 & 100.0 & 100.0 \\
		1000 & $\sigma_1$ & 100.0 & 100.0 & 100.0 & 100.0 & 100.0  & 100.0 & 100.0 & 100.0 & 100.0 & 100.0 \\
		1000 & $\sigma_2$ & 100.0 & 100.0 & 100.0 & 100.0 & 100.0  & 100.0 & 100.0 & 100.0 & 100.0 & 100.0 \\
		1000 & $\sigma_3$ & 100.0 & 100.0 & 100.0 & 100.0 & 100.0  & 100.0 & 100.0 & 100.0 & 100.0 & 100.0 \\
		\midrule
		&&\multicolumn{5}{l}{\textit{Panel C: AR errors $\eps$}} & \multicolumn{5}{|l}{\textit{Panel D: locally stationary errors $\eps$}}\\
		100 & $\sigma_0$ & 73.5 & 100.0 & 100.0 & 100.0 & 98.7  & 90.5 & 100.0 & 100.0 & 100.0 & 100.0 \\
		100 & $\sigma_1$ & 71.1 & 99.8 & 100.0 & 98.2 & 94.1  & 88.1 & 100.0 & 100.0 & 100.0 & 99.1 \\
		100 & $\sigma_2$ & 71.2 & 99.7 & 99.8 & 99.0 & 97.7  & 83.7 & 100.0 & 100.0 & 100.0 & 99.5 \\
		100 & $\sigma_3$ & 69.1 & 99.5 & 99.4 & 99.0 & 96.6  & 80.4 & 99.9 & 99.9 & 99.8 & 98.9 \\
		200 & $\sigma_0$ & 87.2 & 100.0 & 100.0 & 100.0 & 100.0  & 97.8 & 100.0 & 100.0 & 100.0 & 100.0 \\
		200 & $\sigma_1$ & 88.2 & 100.0 & 100.0 & 100.0 & 99.9  & 97.5 & 100.0 & 100.0 & 100.0 & 100.0 \\
		200 & $\sigma_2$ & 85.7 & 100.0 & 100.0 & 100.0 & 100.0  & 92.8 & 100.0 & 100.0 & 100.0 & 100.0 \\
		200 & $\sigma_3$ & 82.5 & 100.0 & 100.0 & 100.0 & 99.9  & 87.9 & 100.0 & 100.0 & 100.0 & 100.0 \\
		500 & $\sigma_0$ & 98.8 & 100.0 & 100.0 & 100.0 & 100.0  & 100.0 & 100.0 & 100.0 & 100.0 & 100.0 \\
		500 & $\sigma_1$ & 99.5 & 100.0 & 100.0 & 100.0 & 100.0  & 100.0 & 100.0 & 100.0 & 100.0 & 100.0 \\
		500 & $\sigma_2$ & 97.4 & 100.0 & 100.0 & 100.0 & 100.0  & 99.7 & 100.0 & 100.0 & 100.0 & 100.0 \\
		500 & $\sigma_3$ & 98.2 & 100.0 & 100.0 & 100.0 & 100.0  & 98.0 & 100.0 & 100.0 & 100.0 & 100.0 \\
		1000 & $\sigma_0$ & 100.0 & 100.0 & 100.0 & 100.0 & 100.0  & 100.0 & 100.0 & 100.0 & 100.0 & 100.0 \\
		1000 & $\sigma_1$ & 100.0 & 100.0 & 100.0 & 100.0 & 100.0  & 100.0 & 100.0 & 100.0 & 100.0 & 100.0 \\
		1000 & $\sigma_2$ & 100.0 & 100.0 & 100.0 & 100.0 & 100.0  & 100.0 & 100.0 & 100.0 & 100.0 & 100.0 \\
		1000 & $\sigma_3$ & 100.0 & 100.0 & 100.0 & 100.0 & 100.0  & 100.0 & 100.0 & 100.0 & 100.0 & 100.0 \\
		\bottomrule
	\end{tabular}
\end{table}

\begin{table}
	\caption{Empirical rejection rates for various choices of $\sigma$ and $\eps$, under the alternative $\mu = \mu_3$ }
	\begin{tabular}{ll|rrrrr|rrrrr}
		\toprule
		$n$ & $\sigma$ & Zhou2013 & CUSUM & \eqref{eq:test_decision_linear} & \eqref{eq:test_decision_poly} & \eqref{eq:test_decision_general} & Zhou2013 & CUSUM & \eqref{eq:test_decision_linear} & \eqref{eq:test_decision_poly} & \eqref{eq:test_decision_general} \\
		\midrule
		&& \multicolumn{5}{l}{\textit{Panel A: i.i.d. errors $\eps$}} & \multicolumn{5}{|l}{\textit{Panel B: MA errors $\eps$}}\\
		100 & $\sigma_0$ & 98.4 & 100.0 & 100.0 & 100.0 & 100.0  & 91.3 & 100.0 & 100.0 & 99.9 & 99.6 \\
		100 & $\sigma_1$ & 98.9 & 100.0 & 100.0 & 100.0 & 100.0  & 90.3 & 100.0 & 100.0 & 99.4 & 97.7 \\
		100 & $\sigma_2$ & 98.8 & 100.0 & 100.0 & 100.0 & 100.0  & 91.9 & 100.0 & 100.0 & 100.0 & 99.7 \\
		100 & $\sigma_3$ & 97.8 & 100.0 & 100.0 & 100.0 & 100.0  & 91.0 & 100.0 & 100.0 & 99.6 & 99.3 \\
		200 & $\sigma_0$ & 100.0 & 100.0 & 100.0 & 100.0 & 100.0  & 98.5 & 100.0 & 100.0 & 100.0 & 100.0 \\
		200 & $\sigma_1$ & 100.0 & 100.0 & 100.0 & 100.0 & 100.0  & 99.0 & 100.0 & 100.0 & 100.0 & 100.0 \\
		200 & $\sigma_2$ & 99.9 & 100.0 & 100.0 & 100.0 & 100.0  & 98.3 & 100.0 & 100.0 & 100.0 & 100.0 \\
		200 & $\sigma_3$ & 99.8 & 100.0 & 100.0 & 100.0 & 100.0  & 97.5 & 100.0 & 100.0 & 100.0 & 100.0 \\
		500 & $\sigma_0$ & 100.0 & 100.0 & 100.0 & 100.0 & 100.0  & 100.0 & 100.0 & 100.0 & 100.0 & 100.0 \\
		500 & $\sigma_1$ & 100.0 & 100.0 & 100.0 & 100.0 & 100.0  & 100.0 & 100.0 & 100.0 & 100.0 & 100.0 \\
		500 & $\sigma_2$ & 100.0 & 100.0 & 100.0 & 100.0 & 100.0  & 100.0 & 100.0 & 100.0 & 100.0 & 100.0 \\
		500 & $\sigma_3$ & 100.0 & 100.0 & 100.0 & 100.0 & 100.0  & 100.0 & 100.0 & 100.0 & 100.0 & 100.0 \\
		1000 & $\sigma_0$ & 100.0 & 100.0 & 100.0 & 100.0 & 100.0  & 100.0 & 100.0 & 100.0 & 100.0 & 100.0 \\
		1000 & $\sigma_1$ & 100.0 & 100.0 & 100.0 & 100.0 & 100.0  & 100.0 & 100.0 & 100.0 & 100.0 & 100.0 \\
		1000 & $\sigma_2$ & 100.0 & 100.0 & 100.0 & 100.0 & 100.0  & 100.0 & 100.0 & 100.0 & 100.0 & 100.0 \\
		1000 & $\sigma_3$ & 100.0 & 100.0 & 100.0 & 100.0 & 100.0  & 100.0 & 100.0 & 100.0 & 100.0 & 100.0 \\
		\midrule
		&&\multicolumn{5}{l}{\textit{Panel C: AR errors $\eps$}} & \multicolumn{5}{|l}{\textit{Panel D: locally stationary errors $\eps$}}\\
		100 & $\sigma_0$ & 80.5 & 100.0 & 100.0 & 100.0 & 99.7  & 95.6 & 100.0 & 100.0 & 100.0 & 100.0 \\
		100 & $\sigma_1$ & 83.6 & 100.0 & 100.0 & 99.1 & 98.7  & 93.6 & 100.0 & 100.0 & 100.0 & 99.9 \\
		100 & $\sigma_2$ & 80.5 & 100.0 & 100.0 & 99.9 & 99.6  & 91.3 & 100.0 & 100.0 & 100.0 & 100.0 \\
		100 & $\sigma_3$ & 81.1 & 100.0 & 100.0 & 99.8 & 99.3  & 84.0 & 100.0 & 100.0 & 99.9 & 99.7 \\
		200 & $\sigma_0$ & 92.7 & 100.0 & 100.0 & 100.0 & 100.0  & 99.8 & 100.0 & 100.0 & 100.0 & 100.0 \\
		200 & $\sigma_1$ & 95.8 & 100.0 & 100.0 & 100.0 & 100.0  & 99.3 & 100.0 & 100.0 & 100.0 & 100.0 \\
		200 & $\sigma_2$ & 91.5 & 100.0 & 100.0 & 100.0 & 100.0  & 97.8 & 100.0 & 100.0 & 100.0 & 100.0 \\
		200 & $\sigma_3$ & 91.4 & 100.0 & 100.0 & 100.0 & 100.0  & 94.2 & 100.0 & 100.0 & 100.0 & 100.0 \\
		500 & $\sigma_0$ & 100.0 & 100.0 & 100.0 & 100.0 & 100.0  & 100.0 & 100.0 & 100.0 & 100.0 & 100.0 \\
		500 & $\sigma_1$ & 100.0 & 100.0 & 100.0 & 100.0 & 100.0  & 100.0 & 100.0 & 100.0 & 100.0 & 100.0 \\
		500 & $\sigma_2$ & 100.0 & 100.0 & 100.0 & 100.0 & 100.0  & 100.0 & 100.0 & 100.0 & 100.0 & 100.0 \\
		500 & $\sigma_3$ & 100.0 & 100.0 & 100.0 & 100.0 & 100.0  & 100.0 & 100.0 & 100.0 & 100.0 & 100.0 \\
		1000 & $\sigma_0$ & 100.0 & 100.0 & 100.0 & 100.0 & 100.0  & 100.0 & 100.0 & 100.0 & 100.0 & 100.0 \\
		1000 & $\sigma_1$ & 100.0 & 100.0 & 100.0 & 100.0 & 100.0  & 100.0 & 100.0 & 100.0 & 100.0 & 100.0 \\
		1000 & $\sigma_2$ & 100.0 & 100.0 & 100.0 & 100.0 & 100.0  & 100.0 & 100.0 & 100.0 & 100.0 & 100.0 \\
		1000 & $\sigma_3$ & 100.0 & 100.0 & 100.0 & 100.0 & 100.0  & 100.0 & 100.0 & 100.0 & 100.0 & 100.0 \\
		\bottomrule
	\end{tabular}
\end{table}

	\end{appendices}

\end{document}